\documentclass{svjour3}
\smartqed  

\usepackage[T1]{fontenc}
\usepackage[latin1]{inputenc}

\journalname{Calcolo}

\usepackage[dvips]{graphicx}
\usepackage{amsmath}
\usepackage{amsfonts}
\usepackage{amssymb}
\usepackage{latexsym}
\usepackage{float}

\usepackage{color}

\def\Mc{{\mathbb M}}

\def\Na{{\mathsf N}}
\def\Ma{{\mathsf M}}
\def\Ka{{\mathsf K}}
\def\Ra{{\mathsf R}}
\def\Pa{{\mathsf P}}
\def\Ca{{\mathsf C}}
\def\Oa{{\mathsf D}}

\def\Nb{\overline{\mathsf N}}
\def\Mb{\overline{\mathsf M}}
\def\Kb{\overline{\mathsf K}}
\def\Rb{\overline{\mathsf R}}
\def\Pb{\overline{\mathsf P}}

\def\Nh{\widehat{\mathsf N}}

\def\Kh{\widehat{\mathsf K}}
\def\Rh{\widehat{\mathsf R}}
\def\Ph{\widehat{\mathsf P}}
\def\Bh{\widehat{\mathsf B}}

\def\xx{\bar{x}}
\def\yy{\bar{y}}

\def\pq{\mathbf q}

\def\rupl{\lceil}
\def\rupr{\rceil}

\def\G{\Gamma}
\def\O{\Omega}

\def\b{\boldsymbol\beta}

\def\g{\boldsymbol{\gamma}}
\def\d{\boldsymbol{\delta}}
\def\l{\lambda}
\def\o{\omega}

\def\et{\boldsymbol{\eta}}

\def\k{\kappa}

\def\h{\hphantom{1}}

\def\qpt{{\quad\forall}}

\def\div{\mathop{\rm div}\nolimits}
\def\rot{\mathop{\rm rot}\nolimits}

\def\dite{\textbf{div}\:}
\def\curl{\textbf{curl}\:}

\newcommand{\lo}{{L^2(\Omega)}}

\newcommand\hcu{{H_{0}^1(\Omega)}}

\newcommand\hcr{H_{0}(\rot;\Omega)}
\newcommand{\N}{{\mathbb{N}}}

\def\bsi{\mbox{\boldmath$\sigma$}}

\def\bta{\mbox{\boldmath$\tau$}}
\def\bep{\mbox{\boldmath$\varepsilon$}}

\def\CT{{\mathcal T}}
\def\CN{{\mathcal N}}

\def\R{\mathbb R}

\def\E{\mathbb E}

\def\bn{\mathbf n}
\def\bt{\mathbf t}

\def\CC{{\mathcal C}}

\def\CE{{\mathcal E}}
\def\CV{{\mathcal V}}
\def\sfv{{\mathsf v}}
\def\sfe{{\mathsf e}}

\def\C{\mathbb C}

\def\tr{\text{tr}}
\def\In{\mathbf{I}}

\newcommand{\ASSUM}[2]{(\textsf{#1#2})}  

\newtheorem{method}{{\bf Method}}

\begin{document}


\title{Numerical results for mimetic discretization of Reissner-Mindlin plate problems}


\author{Louren\c co Beir\~ao da Veiga \and Carlo Lovadina \and David Mora}

\institute{Louren\c co Beir\~ao da Veiga \at
              Dipartimento di Matematica "F. Enriques", Universit\`a degli Studi di Milano,\\
Via Saldini 50, 20133 Milano, Italy \\
              \email{lourenco.beirao@unimi.it}
\and
Carlo Lovadina \at Dipartimento di Matematica, Universit\`a di Pavia,\\
Via Ferrata 1, I-27100 Pavia, Italy\\
\email{carlo.lovadina@unipv.it}
\and
David Mora \at Departamento de Matem\'atica, Universidad del B\'io-B\'io,
Casilla 5-C, Concepci\'on, Chile and CI$^2$MA, Universidad de Concepci\'on, Concepci\'on, Chile \\
\email{dmora@ubiobio.cl}
}

\date{}

\maketitle

\begin{abstract}
A low-order mimetic finite difference (MFD) method for
Reissner-Mindlin plate problems is considered. Together
with the source problem, the free vibration and the
buckling problems are investigated. Full details about
the scheme implementation are provided, and the numerical
results on several different types of meshes are reported.
\end{abstract}

\section{Introduction}\label{intro}

The Reissner-Mindlin theory is widely used to describe the bending
behavior of an elastic plate loaded by a transverse force. However,
its discretization by means of Galerkin methods is typically not
straightforward. For instance, standard low-order finite element
schemes exhibit a severe lack of convergence whenever the thickness
is too small with respect to the other characteristic dimensions of the
plate. This undesirable phenomenon, known as {\em shear locking},
is nowadays well understood: as the plate thickness tends to zero,
the Reissner--Mindlin model enforces the Kirchhoff constraint,
which is typically too severe for the discrete scheme at hand,
especially if low-order polynomials are employed (see, for
instance, the monograph by Brezzi and Fortin~\cite{BF}).
The root of the shear locking phenomenon is that
{\em the space of discrete functions which satisfy the Kirchhoff
constraint is very small, and does not properly approximate a
generic plate solution}.
The most popular way to overcome the shear locking phenomenon
in Galerkin methods  is to reduce the influence of the shear energy by
considering a (selective) reduced integration of the shear part,
by resorting to a mixed formulation or by introducing a suitable
{\em shear reduction operator}. Indeed, several families of methods
have been rigorously shown to be free from locking and optimally convergent;
let us mention, for instance,
\cite{AF,AurLov2001,BD,BeiraodaVeiga:SINUM:2004,BBF:IJNME:1989,BFS:M3AS:1991,DHLRS,DL},
and \cite{L1,L2,Lyly-Niiranen-Stenberg:M3AS:2006,Peisker-Braess:M2AN:1992,Pitkaranta-Suri:NM:1996}.

In the last years, many mimetic discretizations have been
developed for the discretization of problems in partial
differential equations. The mimetic finite difference (or MFD)
method has been successfully employed for solving problems of
electromagnetism \cite{Lipnikov-Manzini-Brezzi-Buffa:JCP:2011},
gas dinamics \cite{Campbell-Shashkov:01}, linear diffusion (see,
e.g.,
\cite{BeiraodaVeiga-Lipnikov-Manzini:NM:2009,BeiraodaVeiga-Lipnikov-Manzini:SINUM:2011,Berndt-Lipnikov-Moulton-Shashkov:01,Berndt-Lipnikov-Shashkov-Wheeler-Yotov:05,BBL09,Brezzi-Lipnikov-Shashkov:05,Brezzi-Lipnikov-Shashkov-Simoncini:07,mr1956021,Hyman-Shashkov-Steinberg:97,Lipnikov-Shashkov-Yotov:09}
and the references therein), convection-diffusion
\cite{Cangiani-Manzini-Russo:09}, Stokes flow
\cite{BeiraodaVeiga-Gyrya-Lipnikov-Manzini:09} and elasticity
\cite{BeiraodaVeiga:M2AN:2010}. We also mention the development of
a posteriori estimators for linear diffusion in
\cite{BeiraodaVeiga-Manzini:IJNME:2008} and post-processing
technique in \cite{Cangiani-Manzini:08}. Finally, the mimetic
discretization method has been shown to share strong similarities
with the finite volume method in
\cite{Droniou:Eymard-Gallouet-Herbin:2010}.

Recently, a mimetic finite difference (MFD) procedure has been
proposed and theoretically analysed for Reissner-Mindlin plates in~\cite{BM}.
The  method, which can be considered as a MFD version of
the MITC and Dur\'an-Liberman elements, combines the excellent
convergence behaviour of the latter schemes with the great
flexibility in handling the mesh of the former approach.
The aim of this paper is to numerically assess the actual
performance of the MFD method, by considering the source problem,
as well as the free vibration and the plate buckling problems.

A brief outline of the paper is as follows. In Section~\ref{s:model}
we recall the Reissner-Mindlin plate problems, together with the
necessary notations. Section~\ref{mimetic} concerns with a
presentation of the numerical scheme proposed and analysed in~\cite{BM}.
In particular, full details about the method implementation,
that where missing in \cite{BM}, are provided. In Section~\ref{s:numer}
we report the numerical results obtained  using several types of meshes.
We end the paper with some concluding remarks.

\section{The Reissner-Mindlin plate equations}\label{s:model}

Here and thereafter we use the following operator notation for any
tensor field $\bta=(\tau_{ij})$ $i,j=1,2$, any vector field
$\et=(\eta_i)$ $i=1,2$ and any scalar field $v$:
\begin{equation*}
\begin{split}
\div\et:=\partial_1\eta_1+\partial_2\eta_2,
\qquad
\rot\et:=\partial_1\eta_2-\partial_2\eta_1,
\qquad
\tr(\bta):=\tau_{11}+\tau_{22},
\\
\nabla v:=\begin{pmatrix}
\partial_1v \\ \partial_2v
\end{pmatrix},
\qquad
\curl v:=\begin{pmatrix}
\partial_2v \\ -\partial_1v
\end{pmatrix},
\qquad
\dite\bta:=\begin{pmatrix}
\partial_1\tau_{11}+\partial_2\tau_{12}
\\
\partial_1\tau_{21}+\partial_2\tau_{22}
\end{pmatrix}.
\end{split}
\end{equation*}

Consider an elastic plate of thickness $t$ such that $0<t\le {\rm
diam}(\Omega)$, with reference configuration
$\O\times\left(-\frac{t}{2},\frac{t}{2}\right),$ where $\O$ is a
convex polygonal domain of $\R^{2}$ occupied by the midsection of
the plate. The deformation of the plate is described by means of
the Reissner-Mindlin model in terms of the rotations
$\b=(\beta_{1},\beta_{2})$ of the fibers initially normal to the
plate's midsurface, the scaled shear stresses
$\g=(\gamma_1,\gamma_2)$, and the transverse displacement $w$.
Assuming that the plate is clamped on its whole boundary
$\partial\O$, the following strong equations describe the plate's
response to conveniently scaled transversal load $g\in\lo$: find
$(\b,w,\g)$ such that
\begin{equation}
\left\{\begin{array}{lll}
-\dite \C\bep(\b)-\g=\boldsymbol{0}&\text{in}&\O,\\
-\div\g=g&\text{in}&\O,\\
\g=\k t^{-2}(\nabla w-\b)&\text{in}&\O,\\
\b=\boldsymbol{0}, w=0&\text{on}&\partial\O ,
\end{array}\right.
\label{msp1}
\end{equation}
where the tensor of bending moduli is given by:
$$\C\bta:=\frac{\E}{12(1-\nu^{2})}\left((1-\nu)\bta+\nu\tr(\bta)\In\right),$$
with $\E>0$ representing the Young modulus, $0< \nu < 1/2$ being
the Poisson ratio for the material and $\In$ indicating the second
order identity tensor.

Let the $\hcu^{2}$-elliptic bilinear form be given by
{\small
\begin{equation}\label{forma}
a(\b,\et):=\int_\O \C \bep(\b) : \bep(\et) =\frac{\E}{12(1-\nu^2)}\int_{\O}\left[(1-\nu)\bep(\b)
:\bep(\et)+\nu\div\b\div\et\right],
\end{equation}}
with $\bep=(\varepsilon_{ij})_{1\le i,j\le2}$ the standard strain
tensor defined by
$\varepsilon_{ij}(\b):=\frac{1}{2}(\partial_{i}\beta_{j}+\partial_{j}\beta_{i}),
1\le i,j\le2$.

Then, the variational formulation of problem \eqref{msp1} reads:
\begin{problem}
\label{erm2}
Find $(\b,w)\in\hcu^{2}\times\hcu$ such that
\begin{equation*}
a(\b,\et)+\k t^{-2}(\nabla w-\b,\nabla v-\et)_{0,\O}
=(g,v)_{0,\O}\qpt(\et,v)\in\hcu^{2}\times\hcu.
\end{equation*}
\end{problem}
In this expression, $\k:=\E k/2(1+\nu)$ is the shear modulus with
$k$ a correction factor usually taken as $5/6$ for clamped plates.

We will also consider the free vibration and buckling problem for plates.

The free vibration problem of a plate is (see \cite{DR,DHELRS,DHLRS,HH}):
\begin{problem}
\label{vrm2}
Find $\l\in\R$ and $0\ne(\b,w)\in\hcu^{2}\times\hcu$ such that
\begin{equation*}
a(\b,\et)+\k t^{-2}(\nabla w-\b,\nabla v-\et)_{0,\O}
=\l \left[(w,v)_{0,\O}+\frac{t^{2}}{12}(\b,\et)_{0,\O}\right]
\end{equation*}
\end{problem}
for all $(\et,v)\in\hcu^{2}\times\hcu$, where $\l =\rho\o^2/t^2$, with $\rho$ being the density
and $\o$ the angular vibration frequency of the plate
and the corresponding eigenfunctions are the vibration modes.

The buckling problem of a plate is (see \cite{KXWL,LMR}):
\begin{problem}
\label{brm2}
Find $\l^{bp} \in\R$ and $0\ne(\b,w)\in\hcu^{2}\times\hcu$ such that
\begin{equation*}
a(\b,\et)+\k t^{-2}(\nabla w-\b,\nabla v-\et)_{0,\O}
=\l^{bp}(\bsi\nabla w,\nabla v)_{0,\O}\qpt(\et,v)\in\hcu^{2}\times\hcu,
\end{equation*}
\end{problem}
where $\bsi(x,y)\in\R^{2\times 2}$ is a symmetric tensor which corresponds
to a pre-existing stress state in the plate,
$\l^{bp}=\l^{bc}/t^{2}$, with $\l^{bc}$ being the buckling coefficients
of the plate and the corresponding eigenfunctions are the buckling modes.

Accordingly with~\eqref{msp1}, for the problems above the scaled shear stresses
can be computed by $\g=\k t^{-2}(\nabla w-\b)$.

\section{A Mimetic Finite Difference (MFD) discretization}\label{mimetic}

In this section we review the mimetic discretization method for the Reissner-Mindlin plate
bending problem presented in \cite{BM}, and extend it to the free vibration and buckling problems.
Finally, in Section~\ref{implementation} we give the details on the implementation of the method.

\subsection{Notation and assumptions}

Let $\{\CT_{h}\}_{h}$ be a sequence of decompositions of the computational domain $\O$ into
$\CN(\CT_h)$ polygons $E$. We assume that this partition is
conformal, i.e. intersection of two different elements $E_1$ and
$E_2$ is either a few mesh vertices, or a few mesh edges (two
adjacent elements may share more than one edge) or empty. We allow
$\CT_h$ to contain non-convex and degenerate elements. For each
polygon $E$, $|E|$ denotes its area, $h_E$ denotes its diameter
and $h:=\max_{E\in\CT_h}h_E.$

We denote the set of mesh vertices and edges by $\CV_h$ and
$\CE_h$, the set of internal vertices and edges by $\CV_h^0$ and
$\CE_h^0$, the set of vertices and edges of a particular element
$E$ by $\CV_h^E$ and $\CE_h^E$, and the set of boundary vertices
and edges by $\CV_h^\partial$ and $\CE_h^\partial$, respectively.
Moreover, we denote a generic mesh vertex by $\sfv$, a generic
edge by $\sfe$ and its length both by $h_\sfe$ and $|\sfe|$.

A fixed orientation is also set for the mesh $\CT_h$, which
is reflected by a unit normal vector $\bn_\sfe$, $\sfe\in\CE_h$, fixed once for all.
Moreover, $\bt_\sfe$ denotes the tangent vector defined as
the counterclockwise rotation of $\bn_\sfe$ by $\pi/2$.

For every polygon $E$ and edge $\sfe\in \CE_h^E$, we define a unit
normal vector $\bn_E^\sfe$ that points outside of $E$, and by
$\bt_E^\sfe$ the tangent vector as the counterclockwise rotation of
$\bn_E^\sfe$ by $\pi/2$.

The mesh is assumed to satisfy the shape regularity
properties detailed in~\cite{BM}.

We make also the following assumption on the data.
The scalar functions $\E, \nu$ are piecewise
constant with respect to the mesh $\CT_h$. Moreover, there exist
two positive constants $C_\star$ and $C^\star$ such that $C_\star
< \E < C^\star$ on the whole domain.
The above uniformity condition on $\E$ is standard, while the
piecewise constant condition can be interpreted as an
approximation of the data and is introduced only for simplicity.
In the general case, it is sufficient to assume that $\E$ and
$\nu$ are (piecewise) $W^{1,\infty}$ and to introduce an
element-wise averaging in the data of the numerical scheme.

\subsection{Degrees of freedom and interpolation operators}

The discretization of Problems~\ref{erm2}-\ref{brm2} requires to discretize
the scalar field of displacement and the vector fields of
rotations and shears. In order to do so, we introduce the degrees
of freedom for the numerical solution in accordance with the
correspondence
\begin{equation*}
\begin{split}
w,v\in\hcu\to w_h,v_h\in W_h,\\
\b,\et\in\hcu^2\to \b_h,\et_h\in H_h,\\
\g,\d\in\lo^2\to \g_h,\d_h\in \G_h,
\end{split}
\end{equation*}
where $W_h$ represents the linear space of discrete displacement,
$H_h$ indicates the linear space of discrete rotations and $\G_h$
is the linear space of discrete shears.

The discrete space for transverse displacements $W_h$ is defined
as follows: a vector $v_h\in W_h$ consists of a collection of
degrees of freedom
$$
v_h:=\{v^\sfv\}_{\sfv\in\CV_h^0},
$$
one per internal mesh vertex, e.g. to every vertex
$\sfv\in\CV_h^0$, we associate a real number $v^\sfv$. The scalar
$v^\sfv$ represents the nodal value of the underlying discrete
scalar field of displacement. The number of unknowns is equal to
the number of internal vertices.

The discrete space for rotations $H_h$ is defined as follows: a
vector $\et_h \in H_h$ is a collection of degrees of freedom
$$
\et_h = \{\et^\sfv\}_{\sfv\in\CV_h^0},
$$
i.e. we assign a vector $\et^\sfv\in\R^2$ per each vertex
$\sfv\in\CV_h^0$.
The vector $\et^\sfv$ represents the nodal values of the
underlying discrete vector field of rotations. The number of unknowns
is equal to twice the number of internal vertices.

Finally, the space for the discrete shear force $\G_h$ is defined
as follows: to every element $E$ in $\CT_h$ and every edge
$\sfe\in\CE_h^E\cap\CE_h^0$, we associate a number
$\delta_E^\sfe$, i.e.
$$
\d_h = \{ \delta_E^\sfe \}_{E\in\CT_h,\sfe\in\CE_h^E\cap\CE_h^0}.
$$
We make the continuity assumption that for each edge $\sfe$ shared
by two element $E_1$ and $E_2$, we have
$$
\delta_{E_1}^\sfe=-\delta_{E_2}^\sfe .
$$
The scalar $\delta_E^\sfe$ represents the average on edges of the
discrete shears in the tangential direction. The number of
unknowns is equal to the number of internal edges.

We now define the following interpolation operators from the
spaces of smooth enough functions to the discrete spaces  $W_h$,
$H_h$ and $\G_h$, respectively. For every function $v\in
\CC^0(\Omega) \cap\hcu$, we define $v_\text{I}\in W_h$ by
\begin{equation}\label{INT0}
v_\text{I}^\sfv:=v(\sfv)\qpt \sfv\in\CV_h^0.
\end{equation}

For every function
$\et\in[\CC^0(\O)\cap\hcu]^2$, we define $\et_{\bf I}\in H_h$ by
\begin{equation}\label{INT1}
\et_{\bf I}^\sfv:=\et(\sfv)\qpt \sfv\in\CV_h^0.
\end{equation}

For every function $\d\in\hcr\cap L^s(\O)^2$, $s>2$, we define
$\d_{\text{II}}\in \G_h$ by
\begin{equation}\label{INT2}
(\delta_\text{II})_E^\sfe:=\frac{1}{\vert \sfe\vert}\int_{\sfe}\d\cdot\bt_{E}^\sfe\qpt E\in\CT_h\qpt \sfe\in\CE_h^E\cap\CE_h^0.
\end{equation}

For all $E\in\CT_h$ in the sequel we will also make use of local
interpolation operators $v_{\text{I},E}, \et_{{\bf I},E}, \d_{\text{II},E}$, with
values in $W_h|_E, H_h|_E, \G_h|_E$ respectively; such operators are simply
the obvious restriction of the global ones to the element $E$ for
functions which are sufficiently regular on $E$.

\begin{remark}\label{rem:1}
We note that in the present paper we are considering the scheme of
\cite{BM} without the edge bubbles, see Remark 4 of \cite{BM}.
Such version of the method is more efficient in terms of accuracy
vs number of degrees of freedom, while the loss of stability is
seen only on very particular mesh patterns. Indeed, in the
numerical test of Section~\ref{s:numer}, only the first family of
(triangular) meshes suffers from such drawback.
\end{remark}

\subsection{Discrete norms and operators}

We endow the space $W_h$ with the following norm
\begin{equation}\label{Q1}
||v_h||_{W_h}^2:=\sum_{E\in\CT_h}||v_h||_{W_h,E}^2=\sum_{E\in\CT_h}|E|\sum_{\sfe\in\CE_h^E}\left[\frac{1}{|\sfe|}(v^{\sfv_2}-v^{\sfv_1})\right]^2,
\end{equation}
where $\sfv_1$ and $\sfv_2$ are the vertices of $\sfe$.
Although irrelevant in~\eqref{Q1}, in the following we will always
consider that $\sfv_1$ and $\sfv_2$, the vertices of a generic
edge $\sfe$, are oriented in such a way that $\bt_{E}^\sfe$ points from
$\sfv_1$ to $\sfv_2$.

In the space $H_h$, we consider the norm
\begin{equation}
|||\et_h|||_{H_h}^2:=\sum_{E\in\CT_h}
|||\et_h|||_{H_h,E}^2=\sum_{E\in\CT_h}|E|\sum_{\sfe\in\CE_h^E}
\left[ \frac{1}{|\sfe|}||\et^{\sfv_1}-\et^{\sfv_2}|| \right]^2,
\end{equation}
where $\sfv_1$ and $\sfv_2$ are the vertices of the edge $\sfe$, and $||\cdot||$ denotes
the euclidean norm on vectors.

In the space $\G_h$, we consider the following norm
\begin{equation}
||\d_h||_{\G_h}^2:=\sum_{E\in\CT_h}
||\d_h||_{\G_h,E}^2=\sum_{E\in\CT_h}|E|\sum_{\sfe\in\CE_h^E}|\delta_E^\sfe|^2
.
\end{equation}
The norms on $W_h$ and $H_h$ are $H^1(\O)$ type discrete
semi-norms, which become norms due to the boundary conditions on
the spaces, while the norm for $\G_h$ is an $L^2(\O)$ type discrete norm.

In the sequel we will also use the following norm on $H_h$, which
is a $|| \bep (\cdot) ||_{0,\O}$ type discrete norm:
\begin{equation}\label{X:norm}
|| \et_h ||_{H_h}^2:=\sum_{E\in\CT_h} ||\et_h||_{H_h,E}^2 =
\sum_{E\in\CT_h} \min_{c \in\R} |||\et_h - c ([- \bar y, \bar
x])_{{\bf I},E} |||_{H_h,E}^2 ,
\end{equation}
where $(\bar x, \bar y)$ are local cartesian coordinates on $E$
which are null on the barycenter of $E$, so that the function $[-
\bar y, \bar x]$ represents a (linearized) rotation around the
barycenter.

We now introduce the discrete gradient operator $\nabla_h$, defined from the set of
nodal unknowns $W_h$ to the set of edge unknowns $\G_h$ as
follows:
$$\nabla_h:W_h\to\G_h$$
$$(\nabla_h v_h)_E^\sfe:=\frac{1}{\vert\sfe\vert}(v^{\sfv_2}-v^{\sfv_1})\qpt E\in\CT_h,\qpt\sfe\in\CE_h^E\cap\CE_h^0,\qpt v_h\in W_h,$$
where $\sfv_1$ and $\sfv_2$ are the vertices of $\sfe$.

We consider also a reduction operator, defined from the discrete
space of rotations $H_h$ to the set of edge unknowns $\G_h$ as
follows:
$$\Pi_h:H_h\to\G_h$$
$$(\Pi_h\et_h)_E^\sfe:=\frac{1}{2}[\et^{\sfv_1}+\et^{\sfv_2}]\cdot\bt_E^\sfe\qpt E\in\CT_h,\qpt \sfe\in\CE_h^E,\qpt \et_h\in H_h,$$
where $\sfv_1$ and $\sfv_2$ are the vertices of $\sfe$.

\subsection{Scalar products and bilinear forms}

We equip the space $\G_h$ with a suitable scalar product,
defined as follows:
\begin{equation}\label{prodgd}
[\g_h,\d_h]_{\G_h}:=\sum_{E\in\CT_h}[\g_h,\d_h]_{\G_h,E},
\end{equation}
where $[\cdot,\cdot]_{\G_h,E}$ is a discrete scalar product on the
element $E$.

The scalar product must satisfy the following stability and consistency
conditions (see~\cite{BM}).
\begin{description}
\item\ASSUM{S}{1} There exist two positive constants $c_{1}$ and
$c_2$ independent of $h$ such that, for every $\d_h\in\G_h$ and
each $E\in\CT_h$, we have
\begin{equation}
c_1||\d_h||_{\G_h,E}^2\le[\d_h,\d_h]_{\G_h,E}\le c_2||\d_h||_{\G_h,E}^2.
\end{equation}\label{s1}
\item\ASSUM{S}{2} For every element $E$, every scalar linear
function $p_1$ on $E$ and every $\d_h\in\G_h$, we have
\begin{equation}
[(\curl p_1)_{\text{II}},\d_h]_{\G_h,E}=\int_{E} p_1
(\rot_{\G_h}\d_h)_{E}-\sum_{\sfe\in\CE_h^E}\delta_E^\sfe\int_{\sfe}
p_1,
\end{equation}\label{s2}
where the operator $(\rot_{\G_h}\d_h)_E:=\frac{1}{|E|}
\sum_{\sfe\in\CE_h^E}\delta_E^\sfe|\sfe|$.
\end{description}

We denote with $a_h(\cdot,\cdot):H_h\times H_h\to\R$ the
discretization of the bilinear form $a(\cdot,\cdot)$, defined as
follows (see \eqref{forma}):
\begin{equation}\label{X2}
a_h(\b_h,\et_h)=\sum_{E\in\CT_h}a_{h}^E(\b_h,\et_h)\qpt\b_h,\et_h\in
H_h ,
\end{equation}
where $a_h^E(\cdot,\cdot)$ is a symmetric bilinear form on each
element $E$, mimicking
$$
a_{h}^E(\b_h,\et_h) \sim \int_E \C \bep(\widetilde \b_h) :\bep (
\widetilde \et_h) \ .
$$

Similarly to the previous case, we introduce a stability and consistency assumptions for
the local bilinear form $a_h^E(\cdot,\cdot)$.
\begin{description}
\item\ASSUM{S}{1$_a$} there exist two positive constants
$\tilde{c}_{1}$ and $\tilde{c}_2$ independent of $h$ such that,
for every $\et_h\in H_h$ and each $E\in\CT_h$, we have
\begin{equation}\label{s1a}
\tilde{c}_1||\et_h||_{H_h,E}^2\le a_h^E(\et_h,\et_h)\le \tilde{c}_2||\et_h||_{H_h,E}^2.
\end{equation}
\end{description}
\begin{description}
\item\ASSUM{S}{2$_a$} For every element $E$, every linear vector
function ${\bf p_1}$ on $E$, and every $\et_h\in H_h$, it holds
\begin{equation}\label{s2a}
a_h^E(({\bf p_1})_{\bf I},\et_h) =
\sum_{\sfe\in\CE_h^E}\Big[\left(\C\bep({\bf p_1})\bn_E^\sfe\right)\cdot\Big(\frac{|\sfe|}{2}[\et^{\sfv_1}+\et^{\sfv_2}]\Big)\Big].
\end{equation}
\end{description}

\begin{remark}\label{rem:2}
The scalar product and the bilinear form shown in this section can
be built element by element in a simple algebraic way.
The details are shown in Section~\ref{implementation}.
\end{remark}

\subsection{The discrete method}

Finally, we are able to define the mimetic discrete
method for Reissner-Mindlin plates proposed in \cite{BM}. Let the loading term
\begin{equation}\label{load}
(g,v_h)_{h}:=\sum_{E\in\CT_h}\bar{g}|_E\sum_{i=1}^{m_E}v^{\sfv_i}\o_E^i,
\end{equation}
where $\sfv_1,\ldots,\sfv_{m_E}$ are the vertices of $E$,
$\bar{g}|_E:=\frac{1}{|E|}\int_{E}g$, and
$\o_E^1,\ldots,\o_E^{m_E}$ are positive weights such that
$$
\int_E p_1 = \sum_{i=1}^{m_E} p_1(\sfv_i) \o_E^i
$$
for all linear functions $p_1$.
The loading term above is an
approximation of
$$ (g,v_h)_{h}\sim\int_{\O}g\tilde v . $$

Then, the discretization of Problem~\ref{erm2} reads:
\begin{method}\label{disc12}
Given $g\in\lo$, find $(\b_h,w_h)\in H_h\times W_h$ such that
\begin{equation*}
a_h(\b_h,\et_h)+\k t^{-2}[\nabla_h w_h-\Pi_h\b_h,\nabla_h
v_h-\Pi_h\et_h]_{\G_h}=(g,v_h)_{h}
\end{equation*}
\end{method}
for all $(\et_h,v_h)\in H_h\times W_h$.

In order to extend the method to the free vibration problem,
we introduce the following mass bilinear form in $H_h \times W_h$
\begin{equation}\label{mass-bil}
m_h(\b_h, w_h ; \et_h, v_h) = \sum_{E\in\CT_h} m_h^E(\b_h, w_h ; \et_h, v_h)
\end{equation}
for all $\b_h , \et_h \in H_h$ and $w_h,v_h \in W_h$, where the local forms
$$
m_h^E(\b_h, w_h ; \et_h, v_h) =
\sum_{i=1}^{m_E} w^{\sfv_i}  v^{\sfv_i} \o_E^i   +
\frac{t^2}{12} \sum_{i=1}^{m_E} ( \b^{\sfv_i} \cdot  \et^{\sfv_i} ) \o_E^i .
$$

Then, the discretization of Problem~\ref{vrm2} reads:
\begin{method}\label{disc-eig}
Find $\l_h\in\R$ and $(\b_h,w_h)\in H_h\times W_h$ such that
\begin{equation*}
a_h(\b_h,\et_h)+\k t^{-2}[\nabla_h w_h-\Pi_h\b_h,\nabla_h
v_h-\Pi_h\et_h]_{\G_h} = \l_h \: m_h(\b_h, w_h ; \et_h, v_h)
\end{equation*}
\end{method}
for all $(\et_h,v_h)\in H_h\times W_h$.

Finally, in order to discretize the buckling problem we introduce a discrete bilinear form
\begin{equation}\label{stress-bil}
b_h(w_h,v_h)=\sum_{E\in\CT_h} b_{h}^E(w_h,v_h) \qpt w_h, v_h\in
W_h ,
\end{equation}
where $b_h^E(\cdot,\cdot)$ is a symmetric bilinear form on each
element $E$, mimicking
$$
b_{h}^E(w_h,v_h) \sim \int_E (\bsi \nabla \widetilde w_h ) \cdot \nabla
\widetilde v_h \ .
$$
We assume for simplicity that the stress datum $\bsi$ is piecewise constant on the mesh,
a condition that can also be considered as an approximation of a given data.
We require that the local bilinear forms $b_h^E (\cdot,\cdot)$ satisfy the following
stability and consistency conditions.
\begin{description}
\item\ASSUM{S}{1$_b$} There exists a positive constant $\hat{c}$
independent of $h$ such that, for every $v_h\in W_h$ and
each $E\in\CT_h$, we have
\begin{equation}
b_{h}^E(v_h,v_h)\le \hat{c}||v_h||_{W_h,E}^2.
\end{equation}
\item\ASSUM{S}{2$_b$} For every element $E$, every scalar linear
function $p_1$ on $E$, and every $v_h\in W_h$, it holds
\begin{equation}
b_h^E((p_1)_\text{I} , v_h) =\sum_{\sfe\in\CE_h^E} \left(\bsi \nabla p_1 \cdot \bn_E^\sfe \right)  \frac{|\sfe|}{2}[v^{\sfv_1} + v^{\sfv_2}],
\end{equation}
where we recall that $\bsi|_E \in \R^{2 \times 2}$ is constant and symmetric.
\end{description}
Such a condition asserts that the discrete bilinear form is exact when tested on linear functions.
We also remark that for the form $b_h^E(\cdot , \cdot)$ we {\em do not} require any lower bound,
such as the one in~\eqref{s1a}. Indeed, assuming a lower bound condition
for $b_h^E( \cdot , \cdot)$ would be unnatural, since the stress
datum $\bsi$ can be a {\em singular} second-order tensor.
%
%

The discretization of Problem~\ref{brm2} then reads:
\begin{method}\label{disc-stress}
Find $\l_h^{bp} \in\R$ and $(\b_h,w_h)\in H_h\times W_h$ such that
\begin{equation*}
a_h(\b_h,\et_h)+\k t^{-2}[\nabla_h w_h-\Pi_h\b_h,\nabla_h
v_h-\Pi_h\et_h]_{\G_h} =
\l_h^{bp} \: b_h(w_h , v_h)
\end{equation*}
\end{method}
for all $(\et_h,v_h)\in H_h\times W_h$.

\subsection{Implementation of the method}\label{implementation}

In this section we describe explicitly how to build the local bilinear forms appearing in the previous sections.

In what follows $m=m(E) \in \N$ will indicate the number of vertices of the polygon $E$.
We number the vertices in counterclockwise sense as $\sfv_1,..,\sfv_m$ and analogously
for the edges $\sfe_1,..,\sfe_m$, so that $\sfv_j$ and $\sfv_{j+1}$ are the endpoints
of edge $\sfe_j$, $j=1,2,..,m$. Note that here and in the sequel all such indexes are
considered modulus $m$, so that the index $m+1$ is identified with the index $1$.
There are a total of $3m$ local degrees of freedom associated to each element of
the mesh, three for each vertex. We order such local degrees of freedom first with
all rotations and then all deflections, ordered as the vertices
$$
\{ \: \et_E^{\sfv_1}, \et_E^{\sfv_2},...,\et_E^{\sfv_m}, v_E^{\sfv_1}, v_E^{\sfv_2},..., v_E^{\sfv_m} \: \},
$$
where $( \et_E, v_E) \in H_h|_E \times W_h |_E$.

The final local bilinear forms $\Mc=\Mc(E) \in \R^{3m \times 3m}$ associated to each element $E$ will be the sum of two parts
\begin{equation}\label{A-0}
\Mc \ = \ \Mc_1 \ + \ \k t^{-2} \Mc_2 ,
\end{equation}
the first one being associated to the $a_h(\cdot,\cdot)$ term and the second one to the shear
energy term. Once the elemental matrices $\Mc$ are built, the global stiffness matrix is
implemented with a standard assembly procedure as in classical finite elements.

\subsubsection{Matrix for the bilinear form $a_h(\cdot,\cdot)$}\label{ss:1}

We start from the bilinear form $a_h(\cdot,\cdot)$, which is the sum of local bilinear forms
that we express as matrices $\Ma=\Ma(E) \in \R^{2 m \times 2m}$
$$
a_{h}^E(\b_E,\et_E) = \b_E^T \Ma \et_E \qpt E \in\CT_h, \ \forall \b_E,\et_E \in H_h |_E.
$$
The first and main step is to build the matrix $\Ma$. With this purpose
we introduce the matrices $\Na=\Na(E)$ and $\Ra=\Ra(E)$ in $\R^{2 m \times 6}$.
Note again that for ease of notation we do not make explicit the dependence on
the involved matrices from $E$. Let $\pq_1,\ldots,\pq_6$ be the following basis
for the first order vector polynomials (with $2$ components) defined on $E$:
$$
\pq_1 = \begin{pmatrix} 1 \\ 0 \end{pmatrix} , \ \pq_2 = \begin{pmatrix} 0 \\ 1 \end{pmatrix} , \
\pq_3 = \begin{pmatrix}  \yy \\ -\xx \end{pmatrix} , \ \pq_4 = \begin{pmatrix} \yy \\ \xx \end{pmatrix} , \
\pq_5 = \begin{pmatrix} \xx \\ \yy \end{pmatrix} , \ \pq_6 = \begin{pmatrix} \xx \\ -\yy \end{pmatrix}.
$$
Then, the six columns $\Na_1,\ldots,\Na_6$ of $\Na$ are vectors in $\R^{2m}$ defined by the
interpolation of the polynomials $\pq_1,\ldots,\pq_6$ into the space $H_h|_E$ (see \eqref{INT1})
$$
\Na_j = (\pq_j)_{{\bf I}, E}
$$
so that for $1 \le i \le m$ and $1 \le j \le 6$
$$
\begin{pmatrix} \Na_{2i-1,j} \\ \Na_{2i,j} \end{pmatrix}  = \pq_j (\sfv_i) .
$$
The columns of the matrix $\Na$ thus represent the linear polynomials $\pq_j$ written in terms of the degrees of freedom of $H_h|_E$.

The columns $\Ra_j$ of the matrix $\Ra$ are instead defined as the vectors in $\R^{2m}$
associated to the right hand side of the consistency condition \ASSUM{S}{2$_a$},
computed with respect to the polynomials $\pq_j$, $j=1,\ldots,6$. In other words $\Ra_j$
is the unique vector in $\R^{2m}$ such that for all $\et_E \in H_h|_E \equiv \R^{2m}$
$$
(\Ra_j)^T  \et_E =
\sum_{i=1}^m \left(\C\bep({\pq_j})\bn_E^{\sfe_i} \right) \cdot \Big(\frac{|\sfe_i|}{2}[\et_E^{\sfv_i}+\et_E^{\sfv_{i+1}}]\Big)
$$
see equation \ASSUM{S}{2$_a$}. Note that, since $\bep({\pq_j}) = {\bf 0}$ for $j=1,2,3$, the first three columns $\Ra_1,\Ra_2,\Ra_3$ of $\Ra$ have all zero entries.

From the definition of the vectors $\Na_j$ and $\Ra_j$, it is clear that the consistency
condition \ASSUM{S}{2$_a$} translates into the algebraic condition
\begin{equation}\label{A-1}
\Ma \Na_j = \Ra_j \ \ j=1,\ldots,6  \quad \Leftrightarrow \quad \Ma \Na = \Ra .
\end{equation}
We therefore introduce the matrix $\Ka \in \R^{6 \times 6}$ defined by
$$
\Ka = \Na^T \Ra = \Ra^T \Na .
$$
It is easy to check that such matrix is symmetric and semi-positive definite. Moreover, it is of the form
$$
\Ka = \begin{pmatrix} {\bf 0}_{3 \times 3} & {\bf 0}_{3 \times 3} \\ {\bf 0}_{3 \times 3}  & \Ka_\star \end{pmatrix}
$$
with $\Ka_\star$ positive definite.  Therefore is it immediate to compute the pseudo inverse of $\Ka$
$$
\Ka^\dagger = \begin{pmatrix} {\bf 0}_{3 \times 3} & {\bf 0}_{3 \times 3} \\ {\bf 0}_{3 \times 3}  & \Ka_\star^{-1} \end{pmatrix} .
$$
We are now ready to define the local matrix $\Ma$. Let $\Pa$ be a projection on the space orthogonal to the columns of $\Na$
$$
\Pa = {\mathbb I}_{2m\times 2m} - \Na (\Na^T \Na)^{-1} \Na^T
$$
with ${\mathbb I}_{2m \times 2m}$ the identity matrix.
We then set
$$
\Ma = \Ra \Ka^\dagger \Ra^T  \ + \ \alpha \Pa
$$
with $\alpha \in \R$ any positive number, typically scaled as the trace of the first part of the matrix.
Then, it is immediate to check that $\Ma$ satisfies the consistency condition \eqref{A-1}.
Moreover, the positivity up to the kernel is easy to check, while the uniform positivity
represented by the stability condition \ASSUM{S}{1$_a$} can be proved with the techniques
shown in \cite{Brezzi-Lipnikov-Simoncini:05,BeiraodaVeiga-Gyrya-Lipnikov-Manzini:09}.

Finally, note that the matrix $\Ma \in \R^{2m \times 2m}$ is defined only with respect to
the rotation degrees of freedom, since the bilinear form $a_h(\cdot,\cdot)$ is
independent of the deflection variable. When it comes to build the local matrix
$\Mc_1 \in \R^{3m \times 3m}$ appearing in \eqref{A-0} one simply needs to
introduce the restriction matrix ${\mathsf S} \in \R^{3m \times 2m}$
$$
{\mathsf S} = \begin{pmatrix} {\mathbb I}_{2m \times 2m} \\ {\bf 0}_{m \times 2m} \end{pmatrix}
$$
and set
$$
\Mc_1 = {\mathsf S} \Ma  {\mathsf S}^T .
$$

\subsubsection{Matrix for the shear term}

The local matrices for the shear part are obtained as a product of matrices representing the operators
and bilinear forms that appear in the second term of the left hand side of Method~\ref{disc12}.
We therefore start building a matrix $\Mb=\Mb(E) \in \R^{m\times m}$ that represents the local scalar product
$$
[\g_E , \d_E ]_{\G_h,E} = \g_E^T \Mb \d_E \quad \forall \g_E,\d_E \in \G_h|_E .
$$
We order the $m$ degrees of freedom of $\G_h|_E$ as the edges of $E$. The construction
follows the same philosophy as in the previous section and therefore is presented more briefly.
Now, the two columns of the matrix $\Nb \in \R^{m \times 2}$ are defined by
$$
\Nb_j = ( \curl q_j)_{\text{II},E} \quad j=1,2 ,
$$
where the sub-index ${}_{\text{II}}$ represents the interpolation operator shown in
\eqref{INT2} and $q_1,q_2$ denote the following basis of the (zero average) linear polynomials on $E$
$$
q_1 = \xx , \ q_2 = \yy .
$$
Analogously, the matrix $\Rb \in \R^{m \times 2}$ is defined through its columns as the right hand side of \ASSUM{S}{2}
$$
(\Rb_j)^T  \d_E = - \sum_{i=1}^m \delta_E^{\sfe_i} \int_{\sfe_i} q_j
\quad \forall j=1,2, \ \forall \d_E \in \G_h|_E \equiv \R^m ,
$$
where we neglected the $\rot_{\G_h}$ part since $q_1$ and $q_2$ have zero average on $E$.
Again, we need to introduce $\Kb \in \R^{2 \times 2}$ given by
$$
\Kb = \Nb^T \Rb = \Rb^T \Nb
$$
that is easily shown to be positive definite and symmetric. We can therefore finally set
$$
\Mb = \Rb \: (\Kb)^{-1} \Rb^T  \ + \ \alpha \Pb
$$
with $\alpha \in \R^+$ and the projection matrix
$$
\Pb = {\mathbb I}_{m\times m} - \Nb (\Nb^T \Nb)^{-1} \Nb^T .
$$
The consistency condition $\Mb\Nb=\Rb$ follows by construction while the stability can
be derived with the results in \cite{Brezzi-Lipnikov-Simoncini:05}.

The local matrix $\Mc_2$ appearing in \eqref{A-0} can be built combining $\Mb$
with a matrix $\Ca=\Ca(E) \in \R^{m \times 3m}$ representing the $\nabla_h$
and $\Pi_h$ operators that appear in Method~\ref{disc12}. We therefore set
$$
\Ca = \begin{pmatrix} - \Ca_1 & \Ca_2 \end{pmatrix}
$$
with the matrix $\Ca_1 = \Ca_1(E) \in \R^{m \times 2m}$ representing the $\Pi_h$ operator
$$
\Ca_1 = \frac{1}{2} \begin{pmatrix}
(\bt_E^{\sfe_1})^T &  (\bt_E^{\sfe_1})^T & {\bf 0}_{1\times 2} & {\bf 0}_{1\times 2} & ...... &  {\bf 0}_{1\times 2} \\
{\bf 0}_{1\times 2} & (\bt_E^{\sfe_2})^T &  (\bt_E^{\sfe_2})^T & {\bf 0}_{1\times 2} & ...... &  {\bf 0}_{1\times 2} \\
{\bf 0}_{1\times 2} & {\bf 0}_{1\times 2} & (\bt_E^{\sfe_3})^T &  (\bt_E^{\sfe_3})^T & ...... &  {\bf 0}_{1\times 2} \\
\vdots & \vdots & \vdots & \vdots & \vdots & \vdots \\
(\bt_E^{\sfe_m})^T  & {\bf 0}_{1\times 2} & ...... & {\bf 0}_{1\times 2} & {\bf 0}_{1\times 2} & (\bt_E^{\sfe_m})^T
\end{pmatrix} ,
$$
and the matrix $\Ca_2 = \Ca_2(E) \in \R^{m \times m}$ representing the $\nabla_h$ operator
$$
\Ca_2 = \begin{pmatrix}
-|\sfe_1|^{-1} &  |\sfe_1|^{-1} & 0 & 0 & ...... &  0 \\
0 & -|\sfe_2|^{-1} &  |\sfe_2|^{-1} & 0 & ...... &  0 \\
0 & 0 & -|\sfe_3|^{-1} &  |\sfe_3|^{-1} &  ...... &  0 \\
\vdots & \vdots & \vdots & \vdots & \vdots & \vdots \\
-|\sfe_m|^{-1} & 0 & ...... & 0 & 0 & |\sfe_m|^{-1}
\end{pmatrix} .
$$
Finally, the local matrices for the shear part are given by
$$
\Mc_2 = \Ca^T \Mb \Ca .
$$

\subsubsection{Right hand sides}

The \emph{loading term} for the source problem in Method~\ref{disc12}
follows immediately from \eqref{load}. One gets the local right hand vectors ${\bf b}={\bf b}(E) \in \R^{3m}$ defined by
$$
{\bf b}_j = \left\{
\begin{aligned}
& 0                                     & \quad \textrm{ if } j=1,2,..,2m   \\
& \bar{g}|_E \: \omega_E^{(j-2m)}  & \quad \textrm{ if } j=2m+1,2m+2,..,3m ,
\end{aligned}
\right.
$$
that are then assembled as usual into the global load vector.

The \emph{mass matrix} for the free vibration problem in Method~\ref{disc-eig},
associated to the bilinear form \eqref{mass-bil} is built again by a
standard assembly procedure. The local mass matrices
$\Oa = \Oa(E) \in \R^{3m \times 3m}$ associated to the elemental mass bilinear forms
$$
\begin{aligned}
m_h^E(\b_E, w_E ; \et_E , v_E) = (\b_E , w_E)^T \Oa \: (\et_E , v_E)& \quad \forall E \in\CT_h,  \\
\forall & \b_E,\et_E \in H_h |_E,  \ \forall w_E, v_E \in W_h|_E
\end{aligned}
$$
are diagonal and defined by
$$
\Oa_{ii} = \left\{
\begin{aligned}
& t^2 \omega_E^{\rupl i/2 \rupr } / 12    & \quad \textrm{ if } i=1,2,..,2m   \\
& \omega_E^{(i-2m)}                   & \quad \textrm{ if } i=2m+1,2m+2,..,3m
\end{aligned}
\right.
$$
where the symbol $\rupl \ \rupr$ stands for a round up to the nearest integer.

The \emph{stress matrix} for the buckling problem in Method~\ref{disc-stress},
associated to the bilinear form \eqref{stress-bil} is also built as the sum of local matrices
$\Bh=\Bh(E) \in \R^{m \times m}$
$$
b_h^E( w_E, v_E) = w_E^T \Bh \: v_E \qquad \forall w_E, v_E \in W_h|_E .
$$
Note that the symmetric tensor $\bsi \in \R^{2 \times 2}$ that appears
in \ASSUM{S}{2$_b$} can have either rank 2 or rank 1. In order to build the matrix
$\Bh$, we start introducing $\{\hat{q}_1,\hat{q}_2,\hat{q}_3 \}$ a basis for the linear polynomials
on $E$, such that $\hat{q}_1=1$ and $\hat{q}_2,\hat{q}_3$ have zero integral on $E$. Moreover,
if $\textrm{rank}(\bsi) = 1$, we also require that $\nabla\hat{q}_2 \in \textrm{Ker}(\bsi)$.
We then define as usual the auxiliary matrices $\Nh=\Nh(E) \in \R^{m \times 3}$
and $\Rh=\Rh(E) \in \R^{m \times 3}$ through its columns. We set
$$
\Nh_j = (\hat{q}_j)_{\text{I},E} ,\quad j=1,2,3 ,
$$
where the sub-index ${}_{\text{I}}$ denotes the interpolation operator in \eqref{INT0},
and define $\Rh_j$ as the unique vector in $\R^m$ such that
$$
\Rh_j^T  v_E
= \sum_{i=1}^m \left(\bsi \nabla\hat{q}_j \cdot \bn_E^{\sfe_i} \right)  \frac{|\sfe_i|}{2}[v_E^{\sfv_i} + v_E^{\sfv_{i+1}}]
\quad \forall j=1,2,3,\;\forall v_E \in W_h|_E \equiv \R^m
$$
in accordance with \ASSUM{S}{2$_b$}. Note that clearly $\Rh_1$ is null, and that,
if $\textrm{rank}(\bsi)=1$ also $\Rh_2$ is null. One then defines as usual the
semi-positive definite and symmetric matrix $\Kh=\Kh(E) \in \R^{3 \times 3}$
$$
\Kh = \Rh^T \Nh = \Nh^T \Rh .
$$
Since $\Kh$ is block diagonal, with the first block of zeros and the second invertible,
it is immediate to compute the pseudo inverse matrix $\Kh^\dagger$,
in a way similar to the one used for $\Ka$ in Section~\ref{ss:1}. Then, we introduce $\Bh=\Bh(E)\in \R^{m \times m}$
$$
\Bh = \Rh \: (\Kh)^\dagger \Rh^T  \ + \ \alpha \Ph
$$
with $\alpha \in \R$ non negative and the projection matrix
$
\Ph = {\mathbb I}_{m\times m} - \Nh (\Nh^T \Nh)^{-1} \Nh^T .
$
Note that, since no global coercivity conditions are required,
differently from the previous matrices also the choice $\alpha=0$ can be taken.

Finally, note that the matrix $\Bh \in \R^{m \times m}$ is defined
only with respect to the deflection degrees of freedom, since the
bilinear form $b_h(\cdot,\cdot)$ is independent of the rotation variable.
The remaining entries in the assembled (right hand side) stress matrix
associated to Method~\ref{disc-stress} can be simply filled with zeros.

\section{Numerical results}\label{s:numer}

The numerical method analyzed has
been implemented in a MATLAB code.

For all the computations we took $\O:=(0,1)^2$, for the Young
modulus we choose: $E=1$.

\begin{figure}[H]
\begin{center}
\begin{minipage}{4cm}
\centering\includegraphics[height=4.5cm, width=4.5cm]{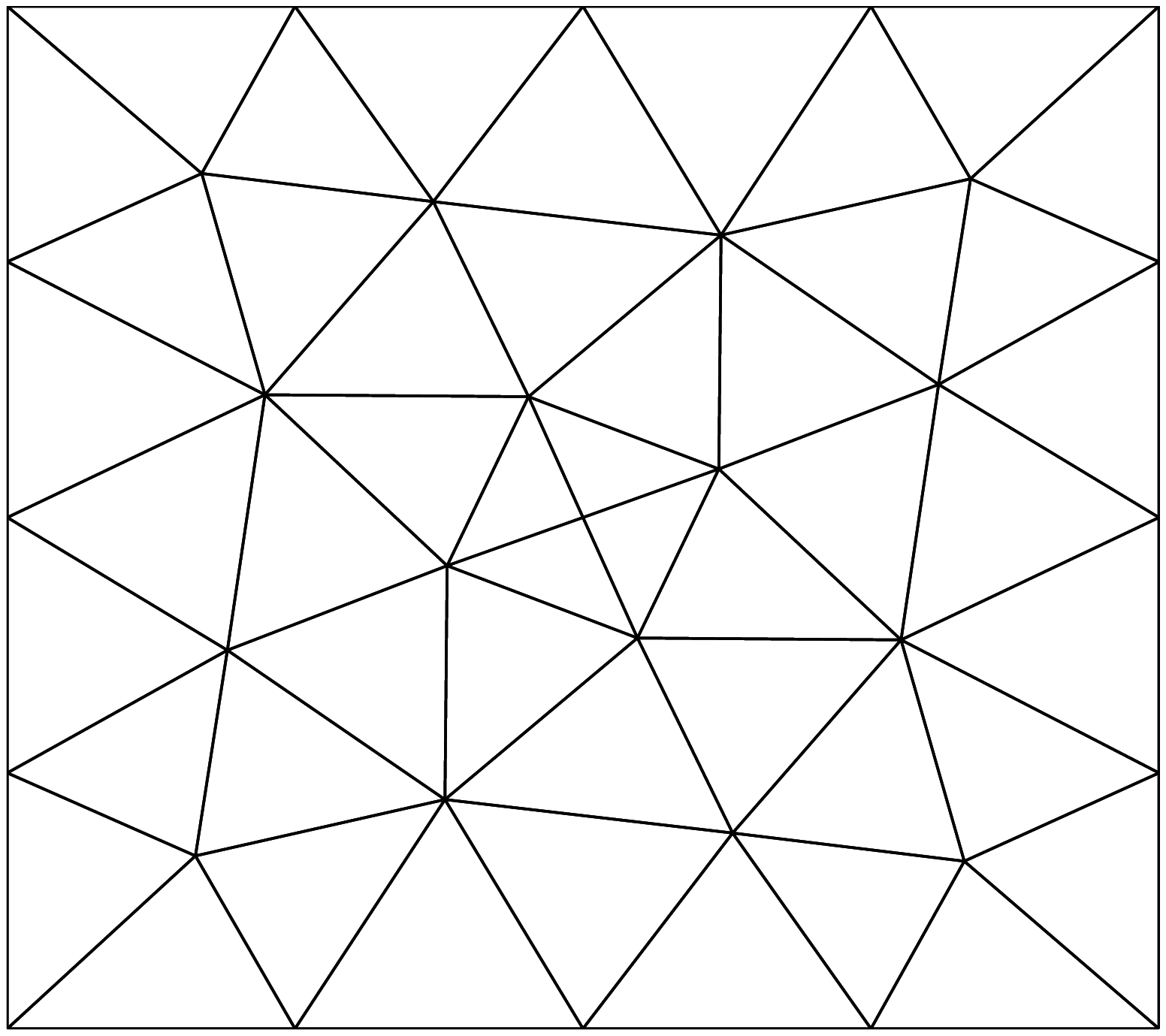}
\centering{$N=4$}
\end{minipage}
\begin{minipage}{3cm}
\centering\includegraphics[height=4.5cm, width=4.5cm]{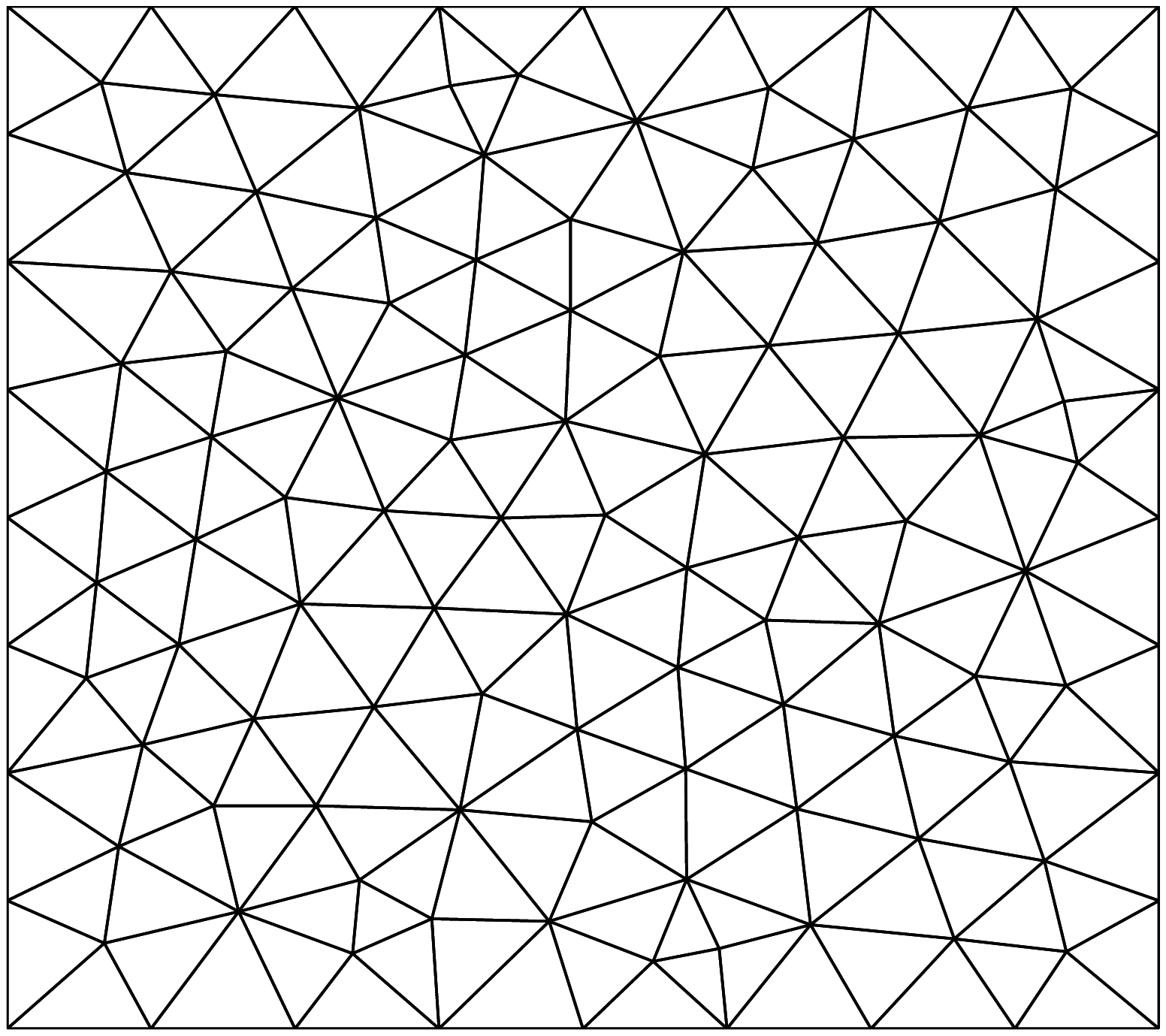}
\centering{$N=8$}
\end{minipage}
\caption{Square plate: meshes $\CT_{h}^{1}$.}
\label{FIG:UNIF_MESH1}
\end{center}
\end{figure}

\begin{figure}[H]
\begin{center}
\begin{minipage}{4cm}
\centering\includegraphics[height=4.5cm, width=4.5cm]{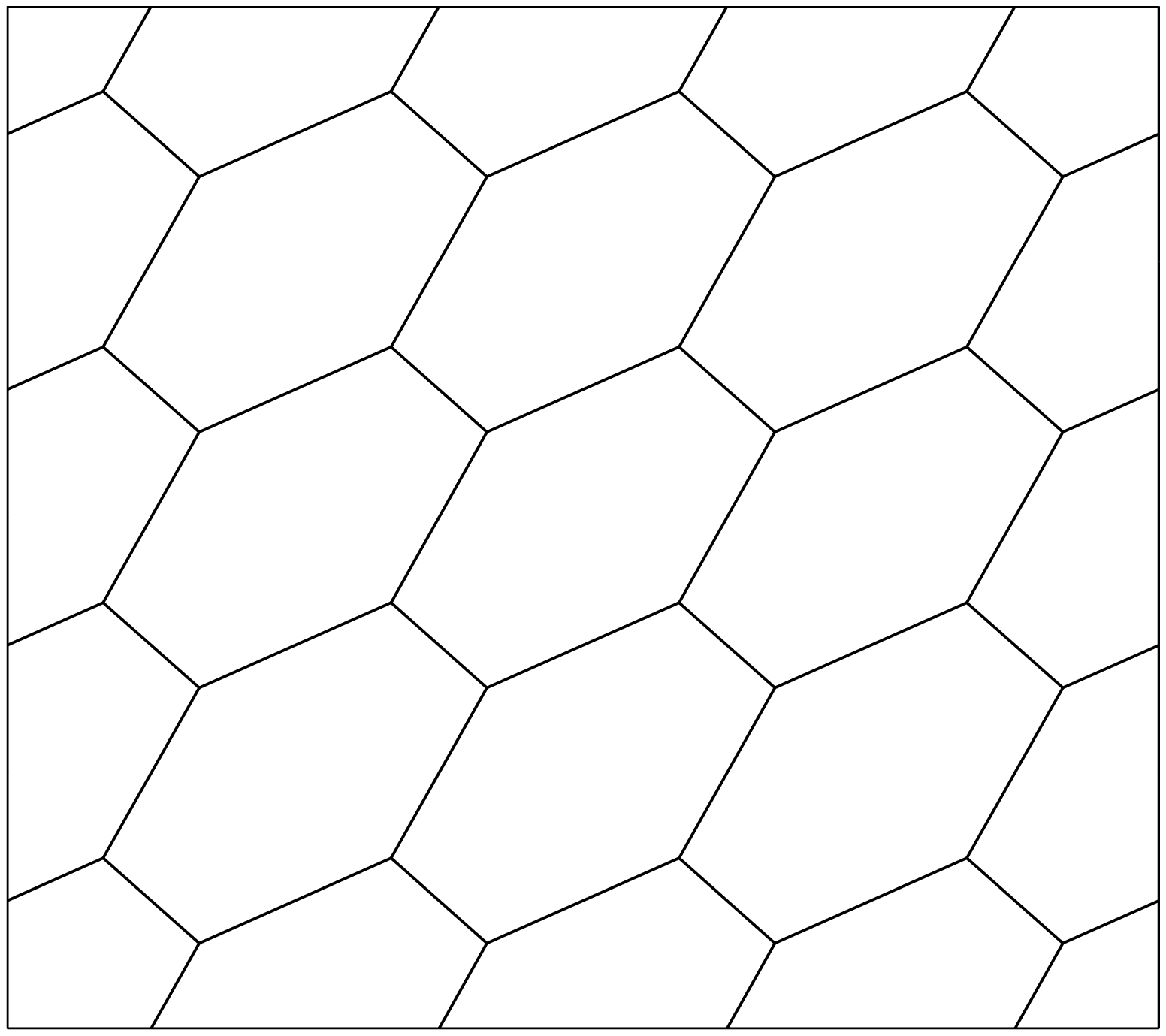}
\centering{$N=4$}
\end{minipage}
\begin{minipage}{3cm}
\centering\includegraphics[height=4.5cm, width=4.5cm]{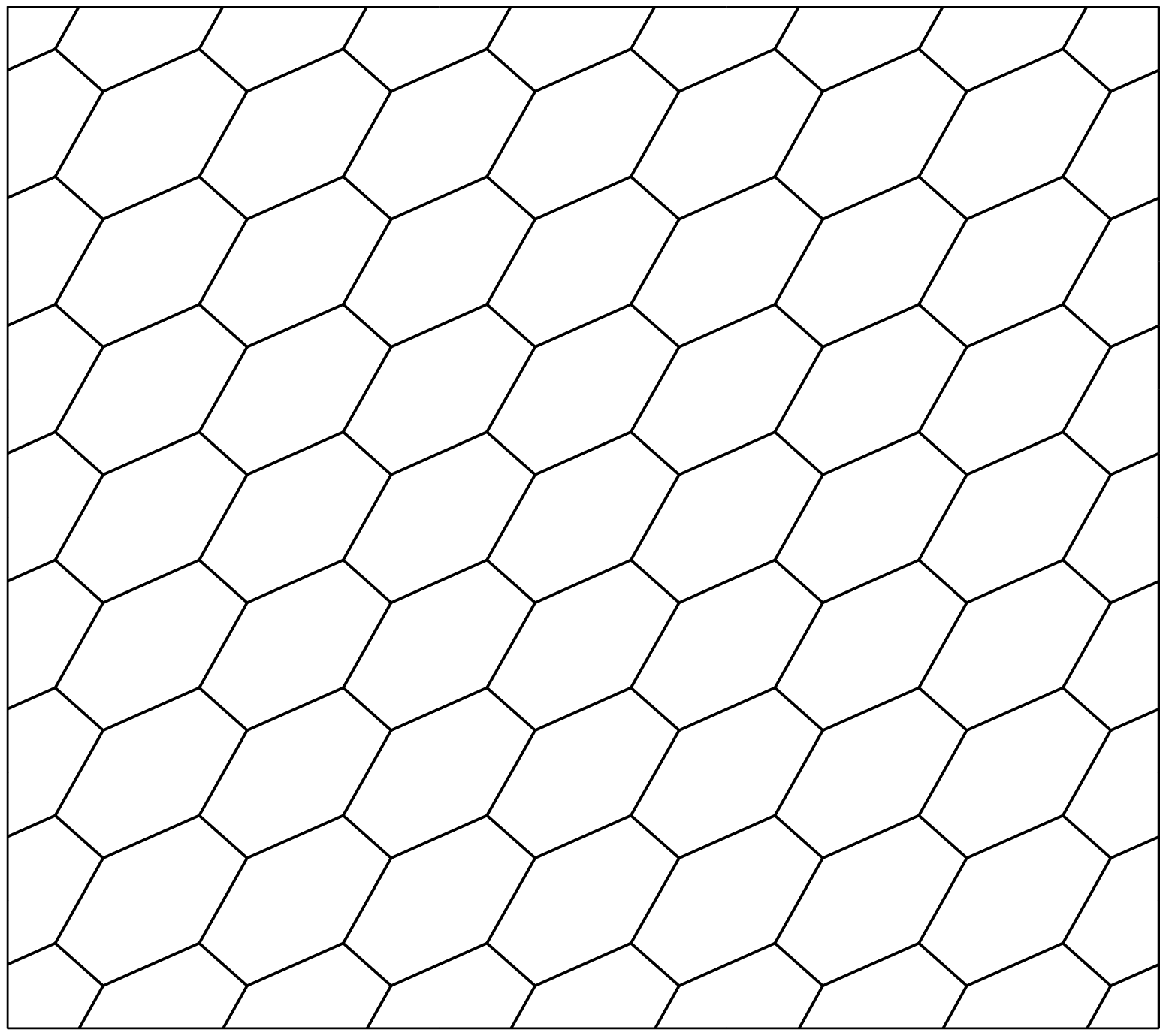}
\centering{$N=8$}
\end{minipage}
\caption{Square plate: meshes $\CT_{h}^{2}$.}
\label{FIG:UNIF_MESH2}
\end{center}
\end{figure}

\begin{figure}[H]
\begin{center}
\begin{minipage}{4cm}
\centering\includegraphics[height=4.5cm, width=4.5cm]{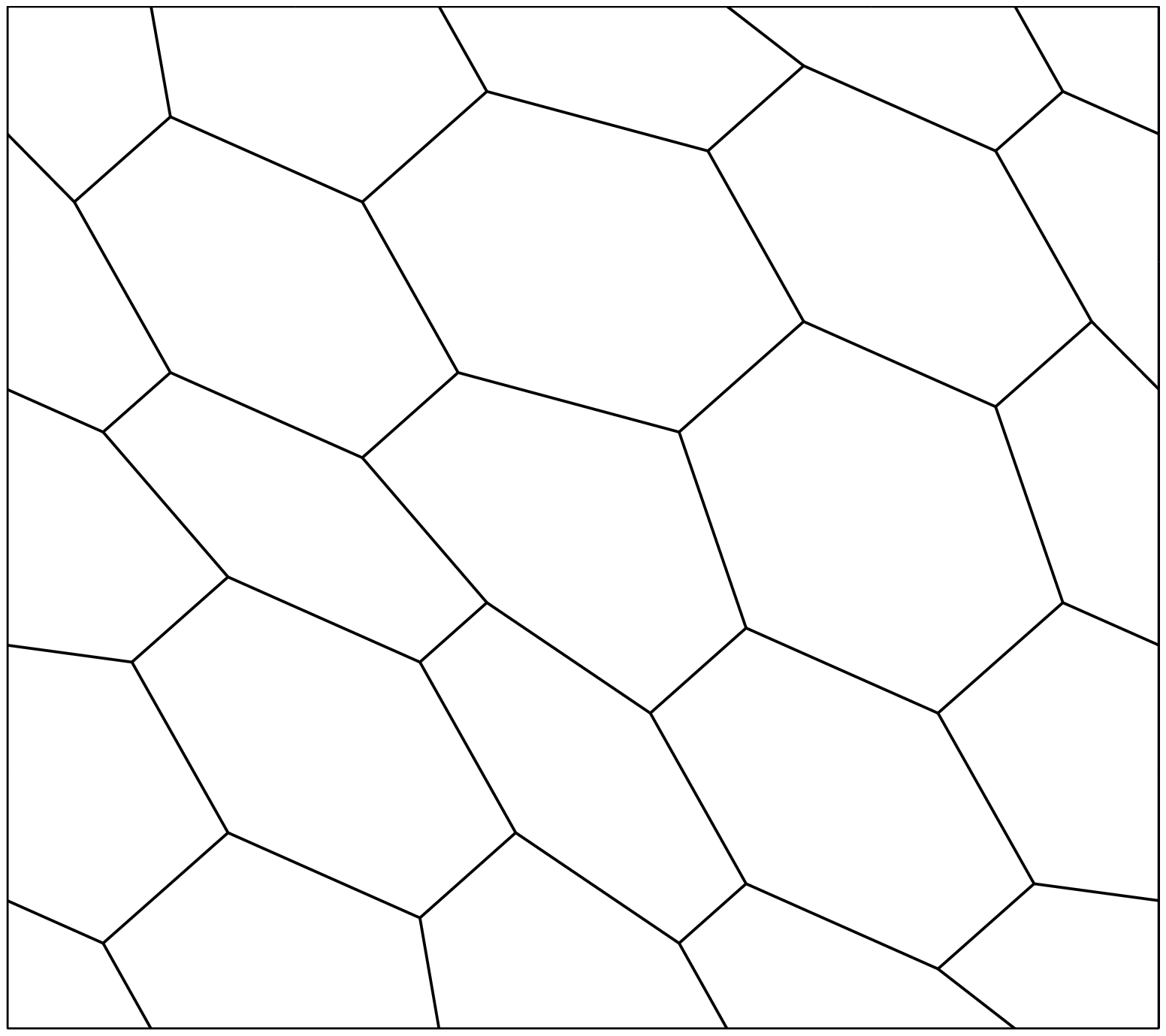}
\centering{$N=4$}
\end{minipage}
\begin{minipage}{3cm}
\centering\includegraphics[height=4.5cm, width=4.5cm]{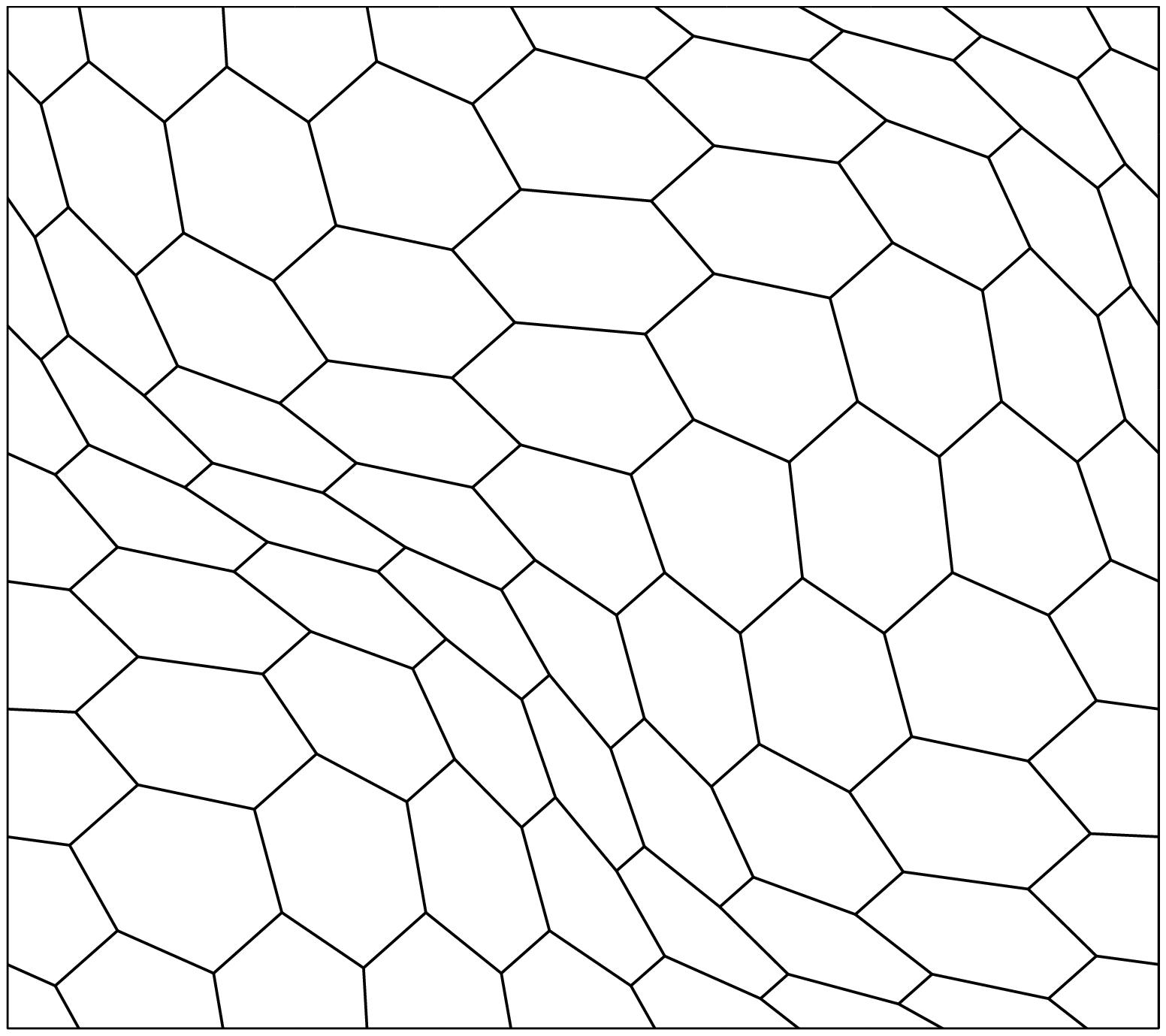}
\centering{$N=8$}
\end{minipage}
\caption{Square plate: meshes $\CT_{h}^{3}$.}
\label{FIG:UNIF_MESH3}
\end{center}
\end{figure}

\begin{figure}[H]
\begin{center}
\begin{minipage}{4cm}
\centering\includegraphics[height=4.5cm, width=4.5cm]{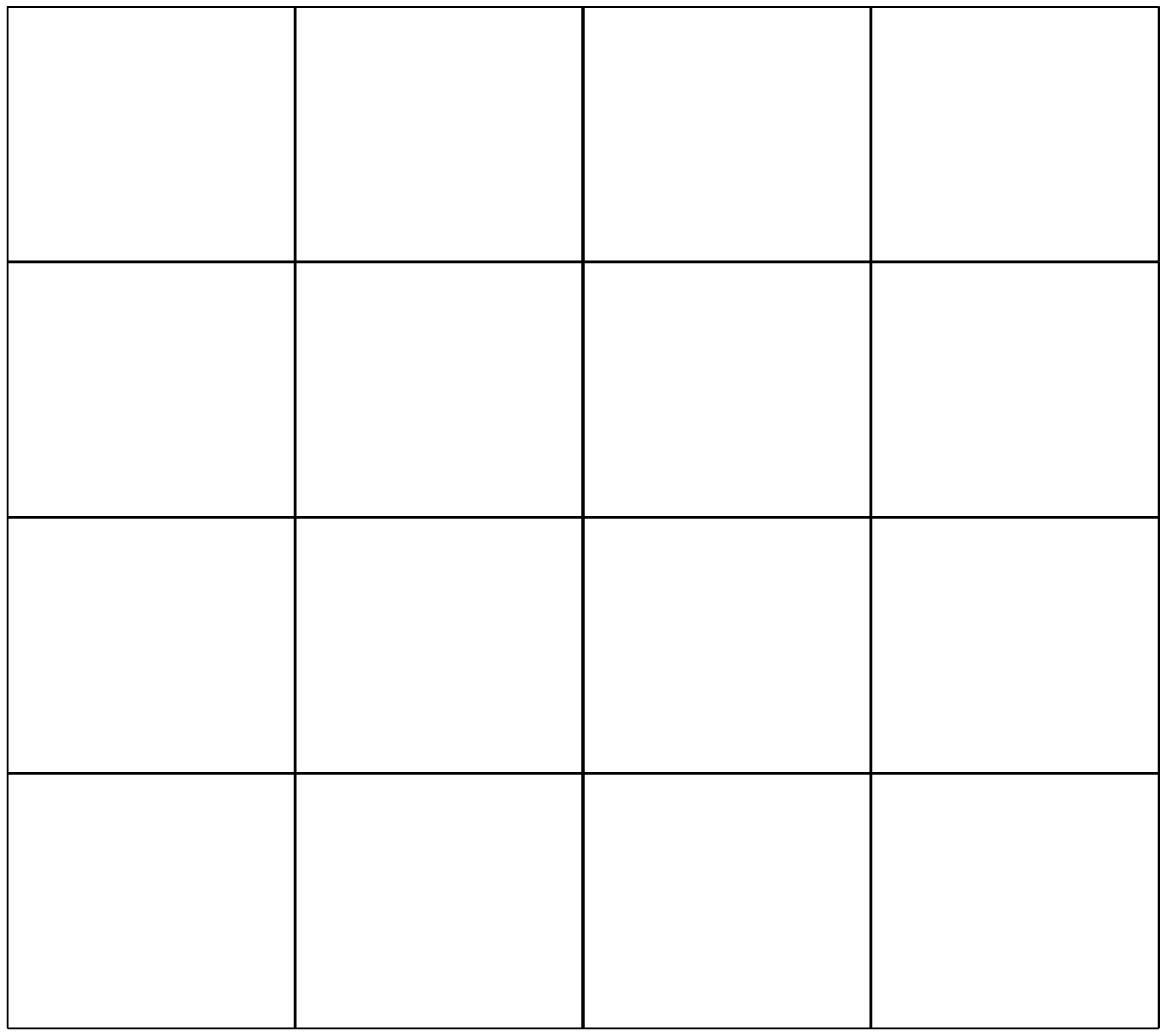}
\centering{$N=4$}
\end{minipage}
\begin{minipage}{3cm}
\centering\includegraphics[height=4.5cm, width=4.5cm]{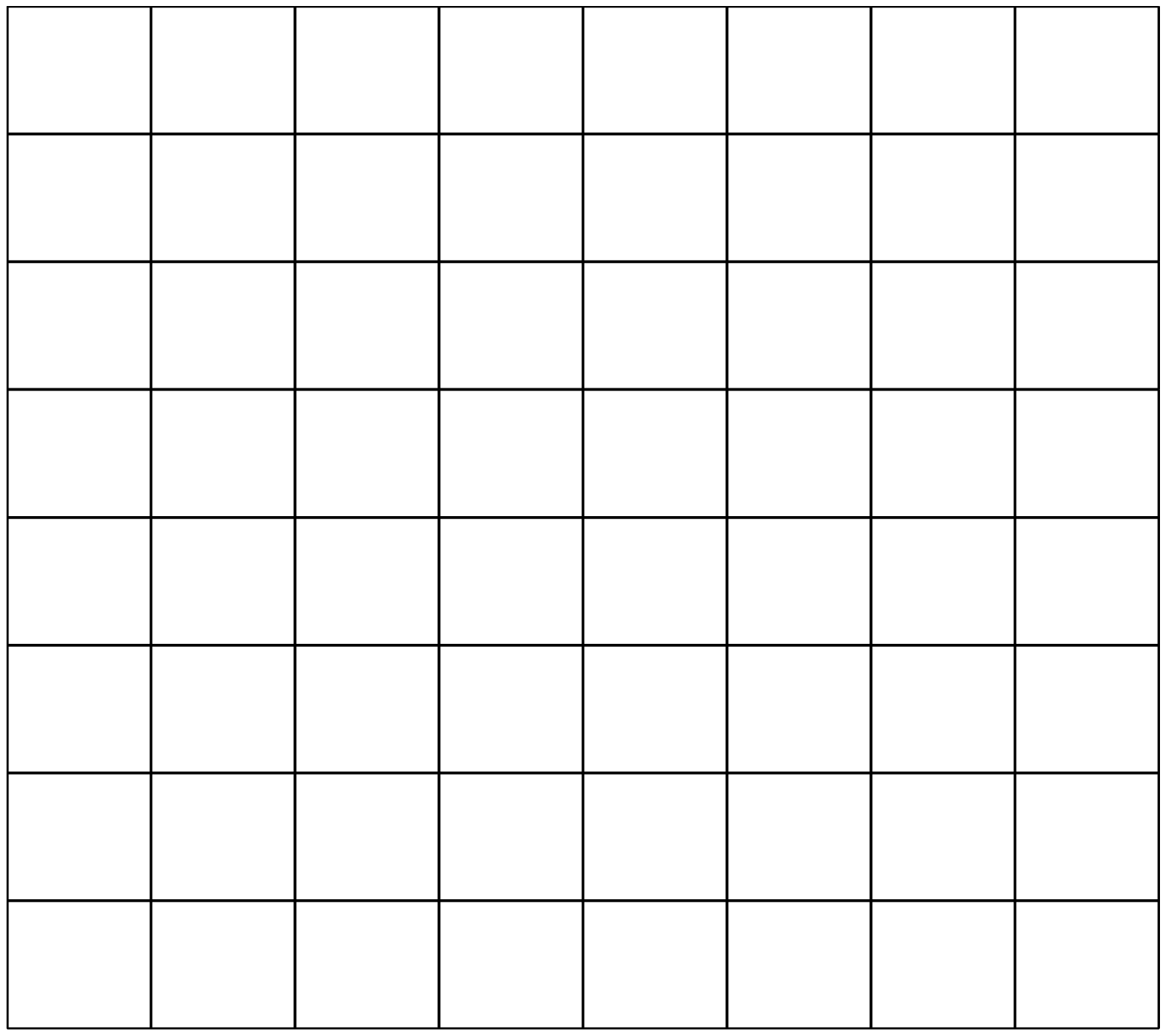}
\centering{$N=8$}
\end{minipage}
\caption{Square plate: meshes $\CT_{h}^{4}$.}
\label{FIG:UNIF_MESH4}
\end{center}
\end{figure}

\begin{figure}[H]
\begin{center}
\begin{minipage}{4cm}
\centering\includegraphics[height=4.5cm, width=4.5cm]{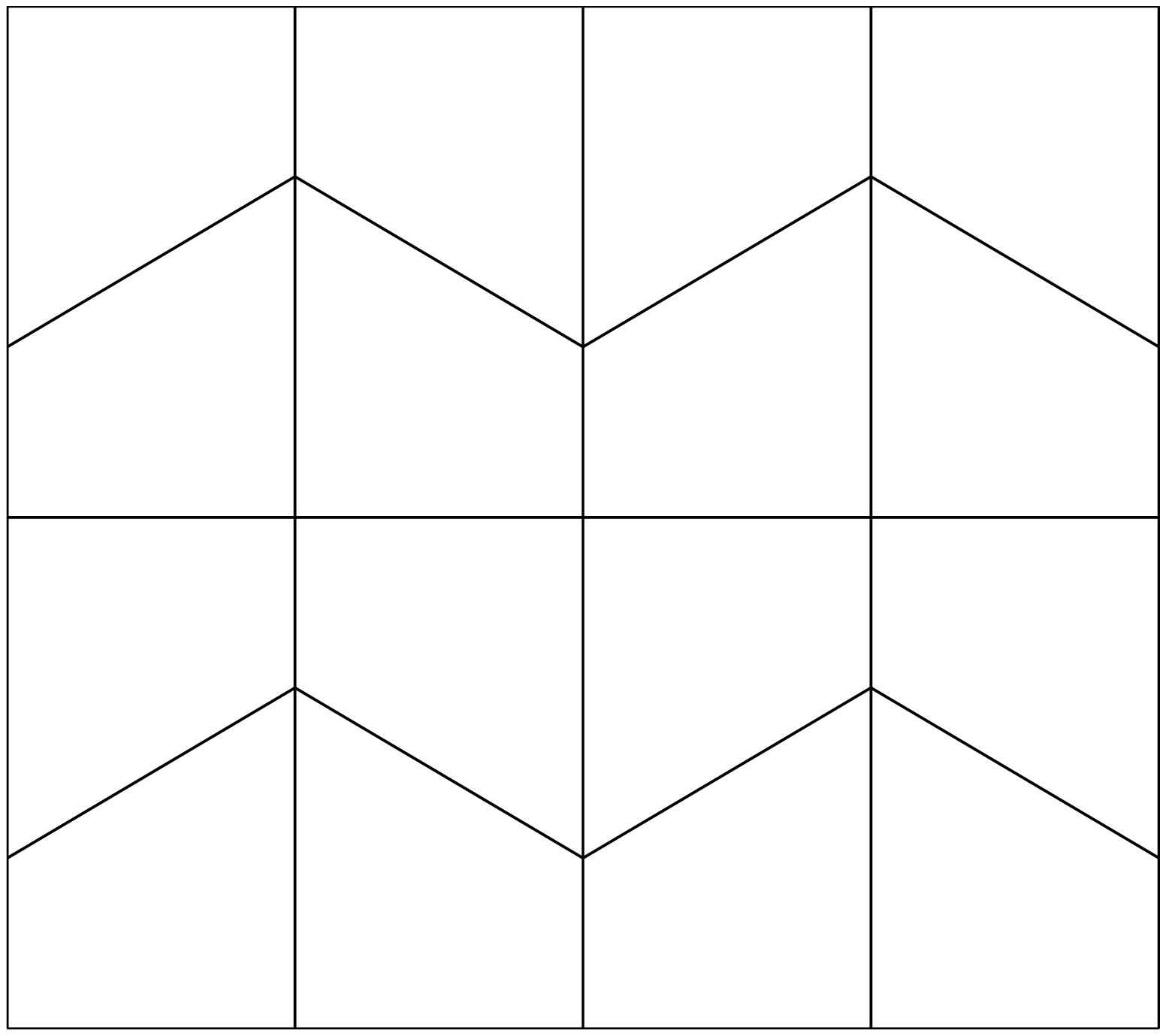}
\centering{$N=4$}
\end{minipage}
\begin{minipage}{3cm}
\centering\includegraphics[height=4.5cm, width=4.5cm]{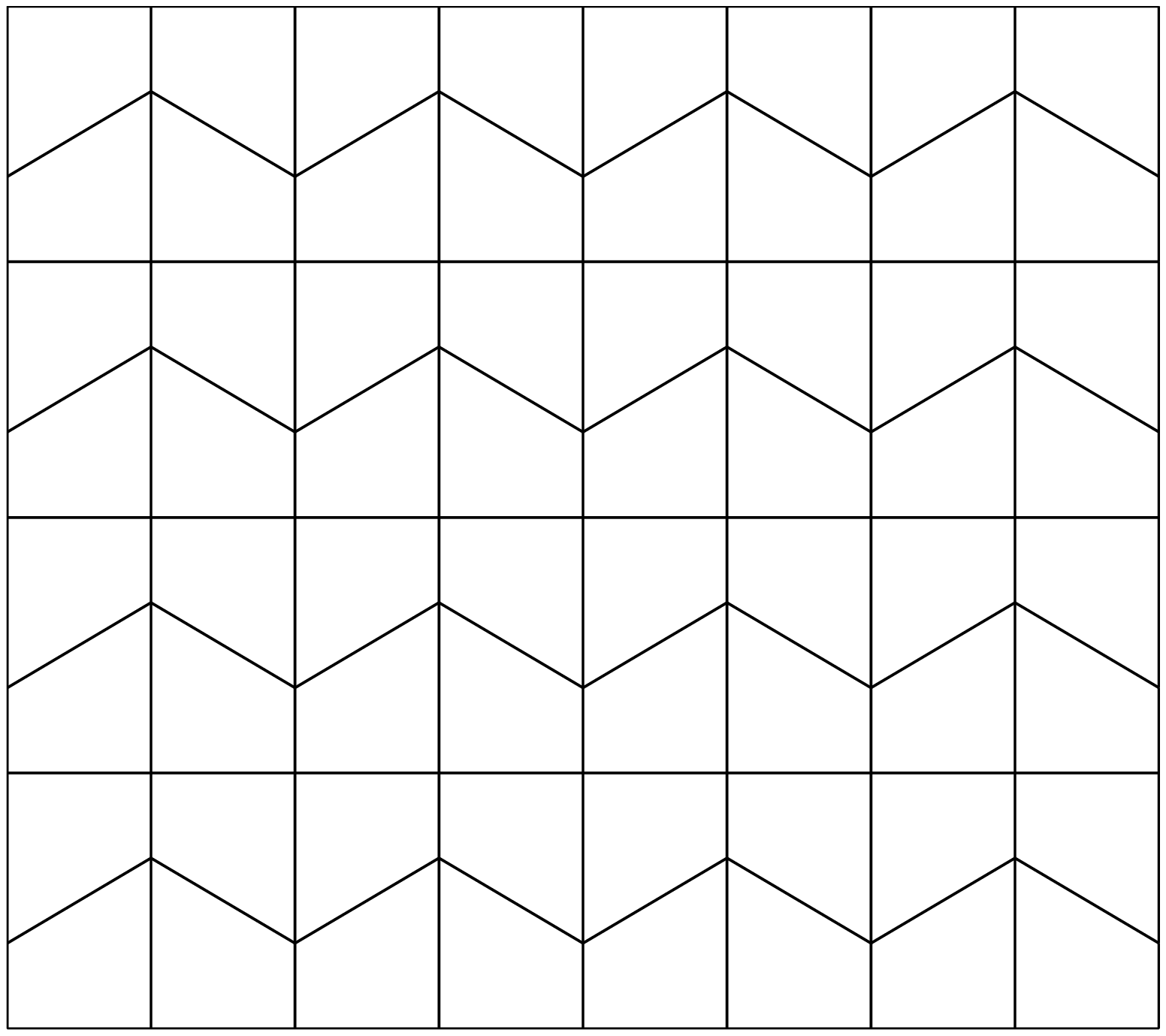}
\centering{$N=8$}
\end{minipage}
\caption{Square plate: meshes $\CT_{h}^{5}$.}
\label{FIG:UNIF_MESH5}
\end{center}
\end{figure}

\begin{figure}[H]
\begin{center}
\begin{minipage}{4cm}
\centering\includegraphics[height=4.5cm, width=4.5cm]{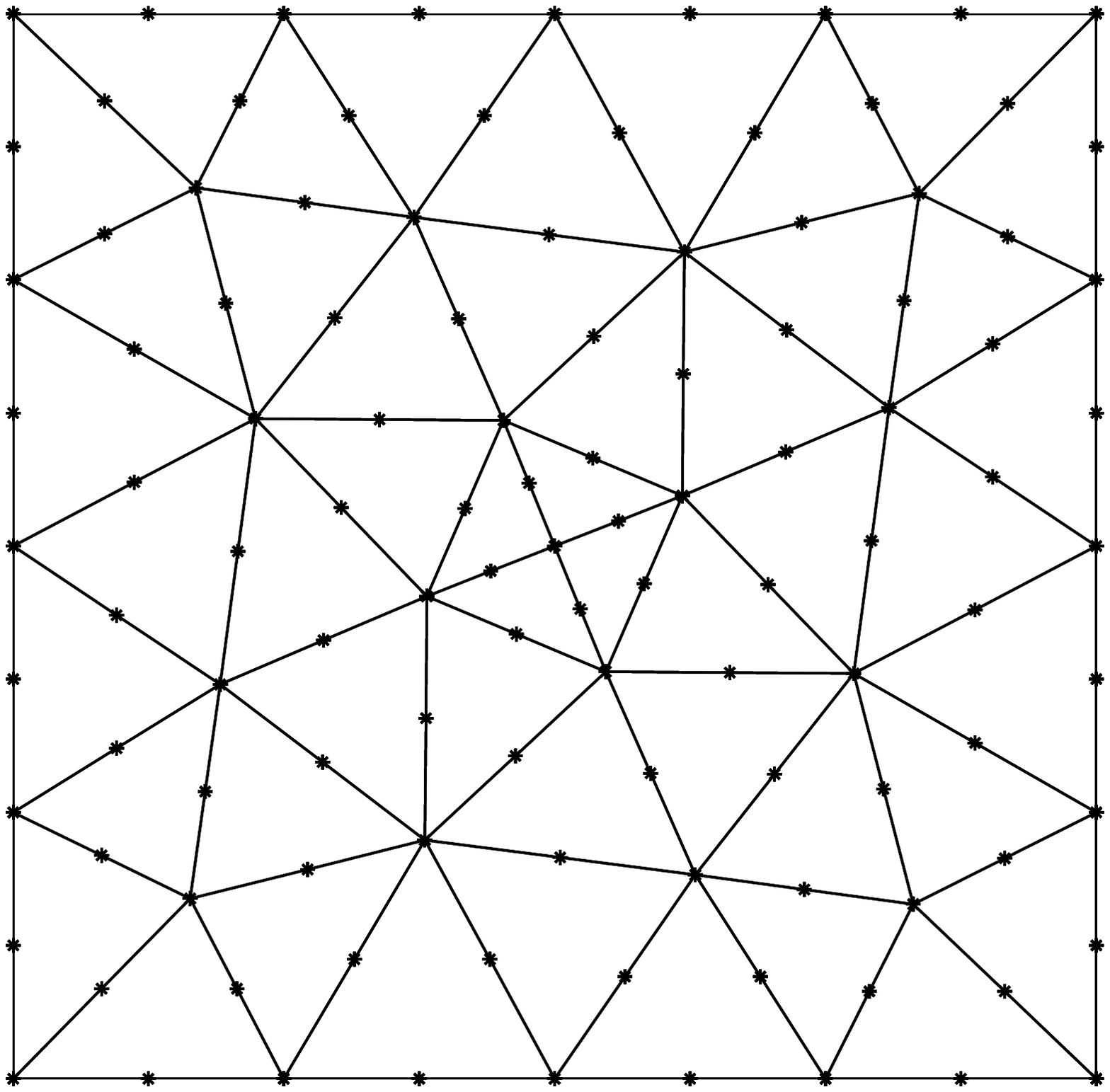}
\centering{$N=4$}
\end{minipage}
\begin{minipage}{3cm}
\centering\includegraphics[height=4.5cm, width=4.5cm]{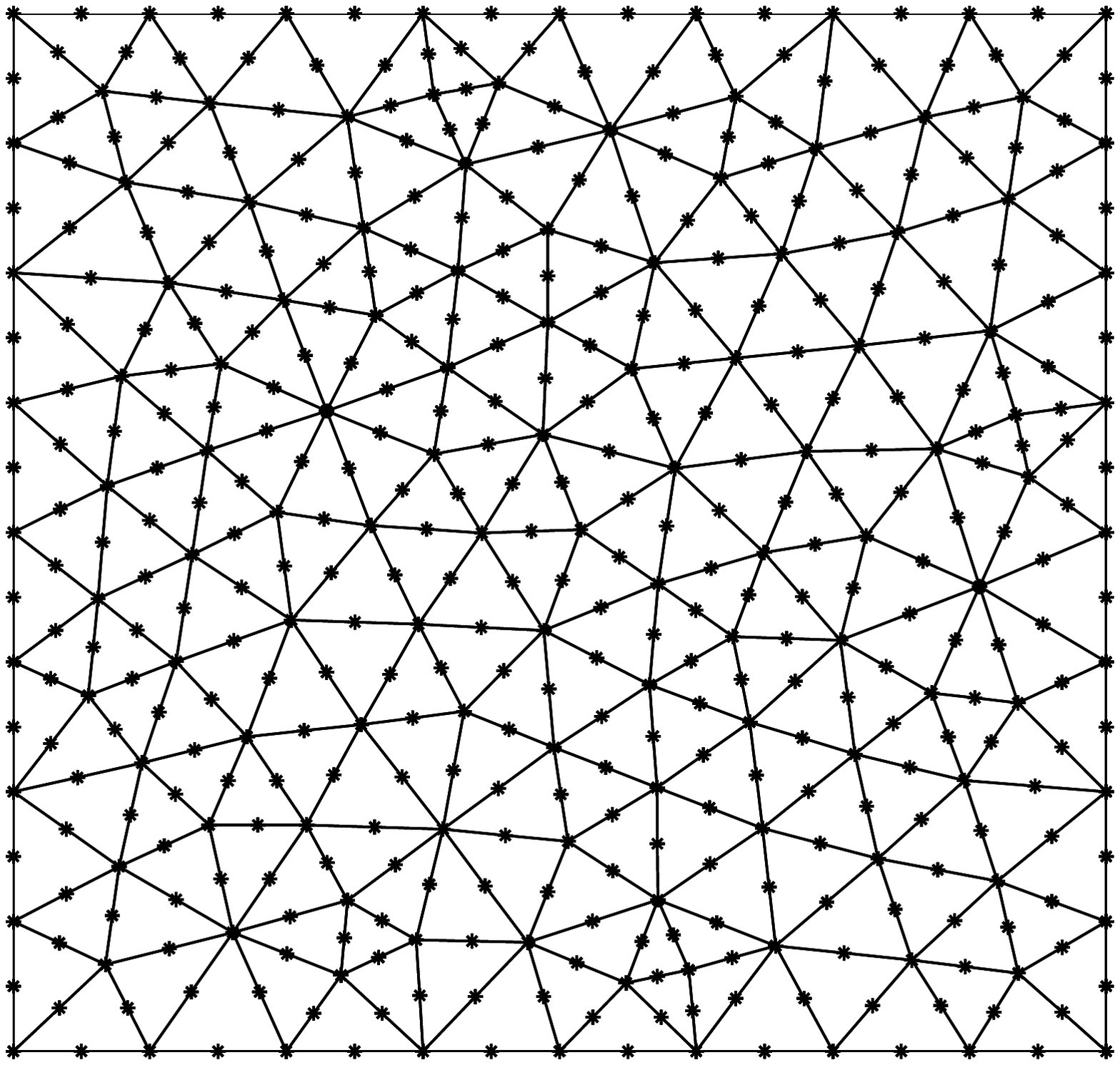}
\centering{$N=8$}
\end{minipage}
\caption{Square plate: meshes $\CT_{h}^{6}$.}
\label{FIG:UNIF_MESH6}
\end{center}
\end{figure}

\begin{figure}[H]
\begin{center}
\begin{minipage}{4cm}
\centering\includegraphics[height=4.5cm, width=4.5cm]{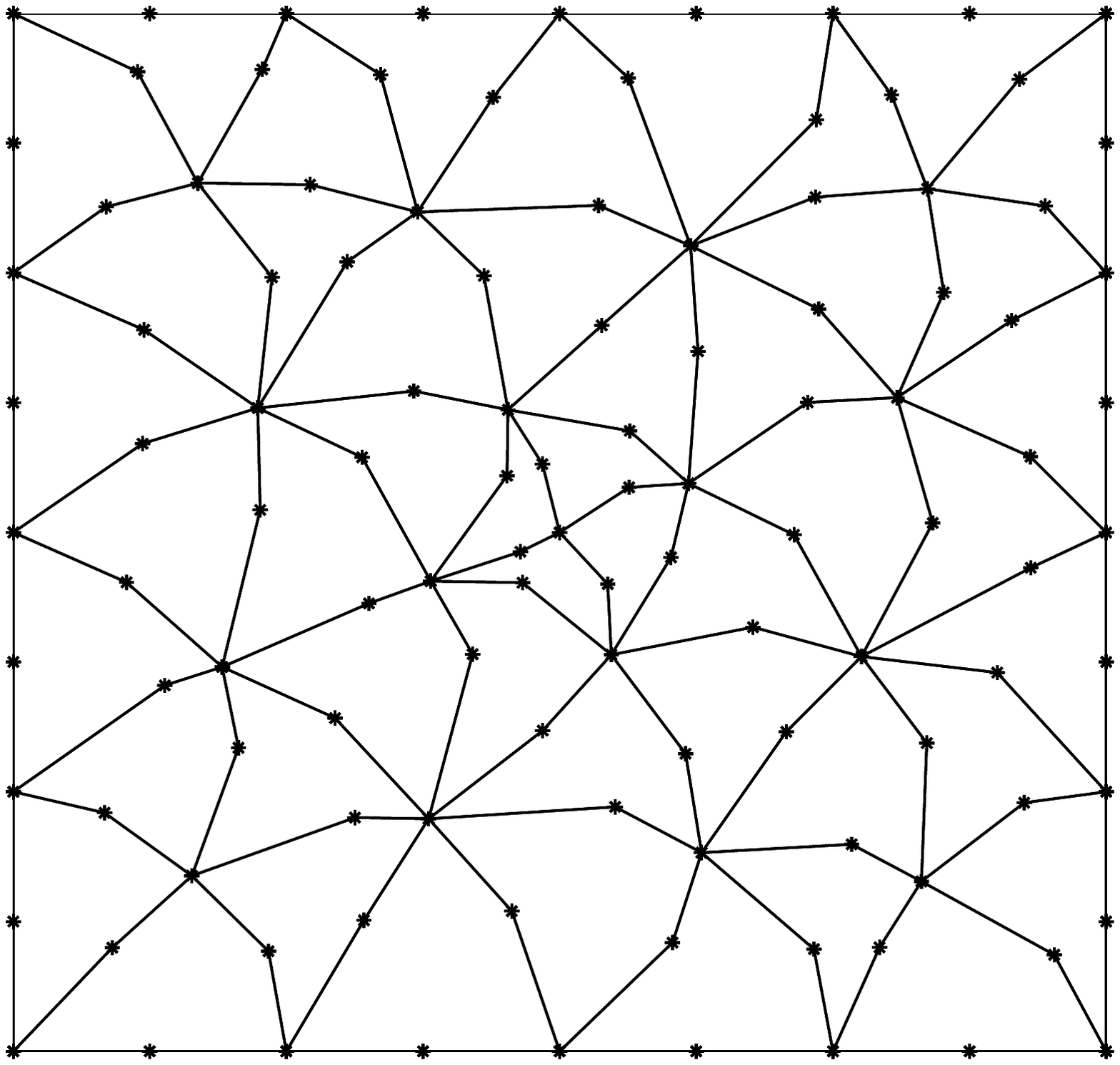}
\centering{$N=4$}
\end{minipage}
\begin{minipage}{3cm}
\centering\includegraphics[height=4.5cm, width=4.5cm]{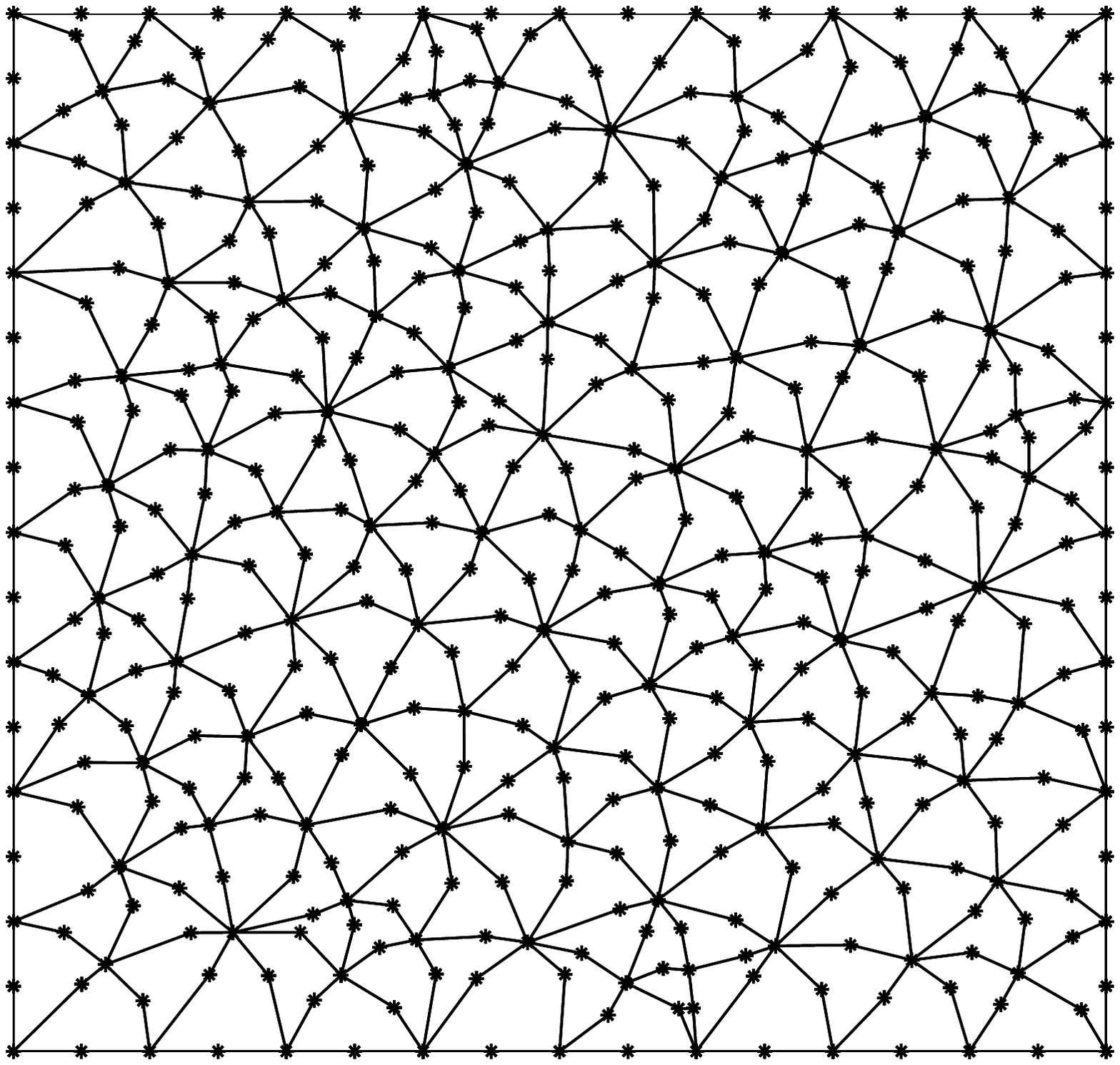}
\centering{$N=8$}
\end{minipage}
\caption{Square plate: meshes $\CT_{h}^{7}$.}
\label{FIG:UNIF_MESH7}
\end{center}
\end{figure}

We have tested the method by using different meshes.
We report the results obtained using the families of meshes
shown in Figure~\ref{FIG:UNIF_MESH1} to Figure~\ref{FIG:UNIF_MESH7}.
\begin{itemize}
\item $\CT_h^1$: Triangular mesh.
\item $\CT_h^2$: Structured hexagonal meshes.
\item $\CT_h^3$: Non-structured hexagonal meshes made of convex hexagons.
\item $\CT_h^4$: Regular subdivisions of the domain in $N\times N$ subsquares.
\item $\CT_h^5$: Trapezoidal meshes which consist of partitions of the domain
into $N\times N$ congruent trapezoids, all similar to the trapezoid with vertices
$(0,0), (\frac{1}{2},0), (\frac{1}{2},\frac{2}{3}),$ and $(0,\frac{1}{3})$.
\item $\CT_h^6$: Regular polygonal meshes built from $\CT_h^1$ considering the middle point
of each edge as a new node on the mesh; note that each element has 6 edges.
\item $\CT_h^7$: Irregular polygonal meshes built from $\CT_h^6$ moving randomly
the middle point of each edge; note that these meshes contain non-convex elements.
\end{itemize}

We have used successive refinements of an initial mesh (see Figure~\ref{FIG:UNIF_MESH1}
to Figure~\ref{FIG:UNIF_MESH7}). The refinement parameter $N$ used to label each mesh
is the number of elements on each edge of the plate.

\subsection{Source problem}

As a test problem we have taken an isotropic and homogeneous plate, clamped
on the whole boundary, for which the analytical solution is explicitly known (see \cite{CLM}).

Choosing the transversal load $g$ as:
{\small
\begin{equation*}
\begin{split}
g(x,y)=\frac{E}{12(1-\nu^2)}\Big[12y(y-1)(5x^{2}-5x+1)
\left(2y^{2}(y-1)^{2}+x(x-1)(5y^{2}-5y+1)\right)\\
+12x(x-1)(5y^{2}-5y+1)
\left(2x^{2}(x-1)^{2}+y(y-1)(5x^{2}-5x+1)\right)\Big],
\end{split}
\end{equation*}
}
the exact solution of the problem is given by:

{\small
\begin{equation*}
\begin{split}
w(x,y)=&\frac{1}{3}x^{3}(x-1)^{3}y^{3}(y-1)^{3}\\
-&\frac{2t^{2}}{5(1-\nu)}\Big[y^{3}(y-1)^{3}x(x-1)(5x^{2}-5x+1)+x^{3}(x-1)^{3}y(y-1)(5y^{2}-5y+1)\Big],
\end{split}
\end{equation*}
\begin{equation*}
\b_{1}(x,y)=y^{3}(y-1)^{3}x^{2}(x-1)^{2}(2x-1),
\end{equation*}
\begin{equation*}
\b_{2}(x,y)=x^{3}(x-1)^{3}y^{2}(y-1)^{2}(2y-1).
\end{equation*}
}

We have used a Poisson ratio $\nu=0.3$ and a correction factor $k=5/6$.

The convergence rates for the transverse displacement $w$
and rotations $\b$ are shown in the following norms:
\begin{equation}\label{error-norm0}
\mathbf{e}(\b)_{0}:=\frac{\max{|\b_{\bf I}-\b_{h}|}}{\max{|\b_{\bf I}|}}
,\quad\mathbf{e}(w)_{0}:=\frac{\max{|w_\text{I}-w_{h}|}}{\max{|w_\text{I}|}} ,
\end{equation}
\begin{equation}\label{error-norm1}
\mathbf{e}(\b)_{1}:=\frac{a_h( \b_{\bf I}-\b_{h} , \b_{\bf I}-\b_{h} )^{1/2}   }{ a_h( \b_{\bf I} , \b_{\bf I})^{1/2} }
,\quad\mathbf{e}(w)_{1}:=\frac{ [ \nabla_h (w_\text{I}-w_{h}) ,  \nabla_h (w_\text{I}-w_{h}) ]_{\G_h}^{1/2}  } { [ \nabla_h w_\text{I} , \nabla_h w_\text{I} ]_{\G_h}^{1/2} }  .
\end{equation}
In~\eqref{error-norm1} the bilinear forms $a_h(\cdot,\cdot)$ and $[ \cdot,\cdot]_{\G_h}$ are exactly the ones defined in~\eqref{X2} and~\eqref{prodgd}, respectively.
We notice that it holds:
\begin{equation*}
\mathbf{e}(\b)_{1} \sim \frac{\Vert \b_{\bf I}-\b_{h}\Vert_{H_h}}{\Vert \b_{\bf I} \Vert_{H_h}}
,\qquad\mathbf{e}(w)_{1} \sim \frac{\Vert w_\text{I}-w_{h}\Vert_{W_h}}{\Vert w_\text{I}\Vert_{W_h}} .
\end{equation*}
Therefore, \eqref{error-norm0} and~\eqref{error-norm1} represent discrete $L^\infty$ and $H^1$ relative errors, respectively.


Also, we define the experimental rates of convergence $(rc)$ for the errors $\mathbf{e}(\b)$ and $ \mathbf{e}(w)$ by
\begin{equation*}
rc(\cdot):=\frac{\log(\mathbf{e}(\cdot)/\mathbf{e'}(\cdot))}{\log(h/h')},
\end{equation*}
where $h$ and $h'$ denote two consecutive meshsizes and $\mathbf{e}$ and $\mathbf{e'}$,
respectively, denote the corresponding errors.

Table~\ref{TAB:L1} shows the
convergence history of the Method~\ref{disc12} applied to our test problem
with four different family of meshes.
Table~\ref{TAB:L3} shows instead an analysis for various thicknesses
in order to assess the locking free nature of the proposed method.

\begin{table}[H]
\caption{Convergence analysis for $t=0.01$.
Errors and experimental rates of convergence for variables $\b$ and $w$.}
\label{TAB:L1}
\begin{center}
{\scriptsize\begin{tabular}{c|c|cc|cc|cc|cc}
\hline
Mesh & $1/h$ & $\mathbf{e}(\b)_{0}$ & $rc(\b)_0$ & $\mathbf{e}(w)_{0}$ & $rc(w)_0$ & $\mathbf{e}(\b)_{1}$ & $rc(\b)_1$&$\mathbf{e}(w)_{1}$ & $rc(w)_1$ \\
\hline
              & 8  &5.018e-2 &--   &2.641e-2 &--   &9.700e-2 &--   &6.480e-2 &-- \\
              & 16 &1.103e-2 &2.19 &9.522e-3 &1.47 &2.967e-2 &1.71 &1.988e-2 &1.70 \\
$\CT_{h}^{2}$ & 32 &2.788e-3 &1.98 &2.761e-3 &1.79 &8.992e-3 &1.72 &5.404e-3 &1.88\\
              & 64 &7.166e-4 &1.96 &7.351e-4 &1.91 &2.886e-3 &1.64 &1.403e-3 &1.95\\
              & 128&1.814e-4 &1.98 &1.891e-4 &1.96 &1.010e-3 &1.51 &3.572e-4 &1.97\\
\hline
              & 8  &7.137e-2 &--   &5.208e-2 &--   &1.325e-1 & --  &8.680e-2 &-- \\
              & 16 &3.505e-2 &1.03 &2.095e-2 &1.31 &5.465e-2 &1.28 &3.303e-2 &1.39\\
$\CT_{h}^{3}$ & 32 &1.131e-2 &1.63 &6.219e-3 &1.75 &1.793e-2 &1.61 &9.714e-3 &1.77\\
              & 64 &3.108e-3 &1.86 &1.634e-3 &1.93 &5.619e-3 &1.67 &2.571e-3 &1.92\\
              & 128&7.991e-4 &1.96 &4.156e-4 &1.98 &1.846e-3 &1.61 &6.571e-4 &1.97 \\
\hline
              & 8  &3.224e-2 &--   &6.519e-2 &--   &4.370e-2 &--   &9.599e-2 &--  \\
              & 16 &8.156e-3 &1.98 &1.605e-2 &2.02 &1.132e-2 &1.95 &2.518e-2 &1.93  \\
$\CT_{h}^{4}$ & 32 &2.051e-3 &1.99 &3.997e-3 &2.00 &2.866e-3 &1.98 &6.365e-3 &1.98  \\
              & 64 &5.138e-4 &1.99 &9.983e-4 &2.00 &7.188e-4 &2.00 &1.595e-3 &2.00  \\
              &128 &1.285e-4 &1.99 &2.496e-4 &2.00 &1.798e-4 &2.00 &3.991e-4 &2.00\\
\hline
              & 8  &7.190e-2 &--   &1.057e-1 &--   &1.949e-1 &--   &1.318e-1 &-- \\
              & 16 &1.677e-2 &2.10 &2.331e-2 &2.18 &1.127e-1 &0.79 &4.339e-2 &1.60 \\
$\CT_{h}^{5}$ & 32 &3.509e-3 &2.26 &5.201e-3 &2.16 &4.213e-2 &1.42 &1.080e-2 &2.01\\
              & 64 &5.942e-4 &2.56 &1.221e-3 &2.09 &1.127e-2 &1.90 &2.380e-3 &2.18\\
              &128 &1.504e-4 &1.98 &2.896e-4 &2.08 &3.802e-3 &1.57 &5.516e-4 &2.11\\
\hline
\end{tabular}}
\end{center}
\end{table}

We observe from Table~\ref{TAB:L1} that a clear rate of convergence
$O(h^2)$ is attained for $\b$ and $w$ for all family of meshes in the
discrete $L^\infty$ norm. Moreover, a rate of convergence
$O(h^{3/2})$ for $\b$ and $O(h^2)$ for $w$ for all family of meshes in the
discrete $H^1$ norm have been obtained.
Actually, the computation of $\b$ using meshes $\CT_{h}^{4}$
seems to be superconvergent.

\begin{table}[H]
\caption{Locking-free analysis for variable $w$ $(\mathbf{e}(w)_{1})$}
\label{TAB:L3}
\begin{center}
{\small\begin{tabular}{c|c|c|c|c|c}
\hline
Mesh & $1/h$ & $t$=1.0e-2 & $t$=1.0e-3 & $t$=1.0e-4 & $t$=1.0e-5 \\
\hline
              & 8   &\h2.040179e-1 &\h9.381597e-1 &\h9.993307e-1 &\h9.999933e-1 \\
$\CT_{h}^{1}$ & 16  &\h2.975046e-2 &\h3.790016e-1 &\h9.773959e-1 &\h9.997674e-1 \\
              & 32  &\h4.320781e-3 &\h4.271249e-2 &\h5.521358e-1 &\h9.902023e-1 \\
              & 64  &\h8.636150e-4 &\h5.651503e-3 &\h8.549716e-2 &\h7.775880e-1 \\
\hline
              & 8   &\h8.680034e-2 &\h8.648290e-2 &\h8.647974e-2 &\h8.647905e-2 \\
$\CT_{h}^{3}$ & 16  &\h3.303805e-2 &\h3.289787e-2 &\h3.289649e-2 &\h3.289827e-2 \\
              & 32  &\h9.714549e-3 &\h9.670538e-3 &\h9.670075e-3 &\h9.671654e-3 \\
              & 64  &\h2.571361e-3 &\h2.558931e-3 &\h2.558797e-3 &\h2.555879e-3 \\
\hline
              & 8   &\h9.598706e-2 &\h9.605054e-2 &\h9.605176e-2 &\h9.605118e-2 \\
$\CT_{h}^{4}$ & 16  &\h2.518265e-2 &\h2.518620e-2 &\h2.518623e-2 &\h2.518663e-2 \\
              & 32  &\h6.364552e-3 &\h6.364126e-3 &\h6.364108e-3 &\h6.364647e-3 \\
              & 64  &\h1.595350e-3 &\h1.595135e-3 &\h1.595115e-3 &\h1.595949e-3 \\
\hline
              & 8   &\h7.471495e-2 &\h7.399270e-2 &\h7.398561e-2 &\h7.398560e-2 \\
$\CT_{h}^{6}$ & 16  &\h2.223311e-2 &\h2.217105e-2 &\h2.217068e-2 &\h2.217124e-2 \\
              & 32  &\h6.155689e-3 &\h6.173926e-3 &\h6.174458e-3 &\h6.176650e-3 \\
              & 64  &\h1.268306e-3 &\h1.308765e-3 &\h1.309856e-3 &\h1.311575e-3 \\
\hline
              & 8   &\h8.776844e-2 &\h8.835775e-2 &\h8.547591e-2 &\h8.870029e-2 \\
$\CT_{h}^{7}$ & 16  &\h3.128646e-2 &\h3.002581e-2 &\h3.088068e-2 &\h2.990004e-2 \\
              & 32  &\h8.776518e-3 &\h8.804238e-3 &\h9.034351e-3 &\h8.770496e-3 \\
              & 64  &\h1.925048e-3 &\h2.025191e-3 &\h2.042556e-3 &\h1.966998e-3 \\
\hline
\end{tabular}}
\end{center}
\end{table}

We observe from Table~\ref{TAB:L3} that our Method~\ref{disc12}
lead to wrong result {\it only} for triangular meshes $\CT_{h}^{1}$,
when the thickness of the plate is small.
For any other family of meshes the method is locking-free.
We also note that adding the middle point of each edge
as a new node on any triangular mesh, the method
is locking free, see row corresponding to $\CT_{h}^{1}$ and $\CT_{h}^{6}$.

\begin{remark}\label{rem:3}
We note that the different behavior among the triangular mesh
$\CT_h^1$ and the remaining grids is not surprising. Indeed,
$\CT_h^1$ resembles a plain $P1$ element, that is known to suffer
from locking in the plate finite element literature, unless edge
bubbles are added to the rotations (see also Remark~\ref{rem:1}).
Moreover, the analysis of \cite{BM} does not apply to the $\CT_h^1$ meshes since
the $P1/P0$ element is not stable for the Stokes problem. 
Instead, meshes $\CT_h^2$, $\CT_h^3$ fall into the hypotheses of
the convergence theorem of \cite{BeiraoLipnikov}, see Remark 4 of
\cite{BM}. Meshes $\CT_h^4$, $\CT_h^5$ have a strong connection
with the MITC4 finite element for plates that is known to be
stable. Finally, meshes $\CT_h^6$, $\CT_h^7$ again fall into the
convergent cases considered in \cite{BeiraoLipnikov} and thus
stable also for Reissner-Mindlin (see \cite{BM}).
\end{remark}

\subsection{Free vibration of plates}

The effectiveness of the MDF method for free vibration analysis are demostrated
by examples with different thickness and different boundary conditions.

We have computed approximations of the free vibration angular frequencies
$\o^{h}=t\sqrt{\frac{\l_{h}}{\rho}}$. In order to compare our results with
those in \cite{DR,DHELRS,DHLRS,HH}, a non-dimensional frequency parameter is defined as:
$$\o_{mn}:=\o_{mn}^{h}L\sqrt{\frac{2(1+\nu)\rho}{E}},$$
here $\o_{mn}^{h}$ are the computed frequencies, where $m$ and $n$ are the numbers
of half-waves in the modal shapes in the $x$ and $y$ directions, respectively.
$L$ is the plate side length.

We have considered a square plate of side length $L=1$ and $\rho=1$ and three different
thickness $t=0.1$, $t=0.01$ and $t$=1.0e-5. We have also considered
three different types of boundary conditions: a clamped plate (denote by CCCC), a simply supported plate
(denote by SSSS), and a plate with a free edge
(with three clamped edges and the fourth free, we denote by CCCF).

In the following numerical tests, we show the results for the four lowest vibration frequencies.
We tested also higher frequencies with similar results.

Tables~\ref{TAB:V1} and \ref{TAB:V2} show the four lowest vibration frequencies computed by
Method~\ref{disc-eig} with successively refined meshes of each type for a clamped plate with thickness
$t=0.1$ and $t=0.01$, respectively. The table includes orders of convegence,
as well as accurate values extrapolated by means of a least-squares fitting.
Furthermore, the last two columns show the results reported in \cite{DR,DHELRS,HH}.
In every case, we have used a Poisson ratio $\nu=0.3$ and a correction
factor $k= 0.8601$. The reported non-dimensional frequencies are independent
of the remaining geometrical and physical parameters, except for the thickness-to-span ratio.

\begin{table}[H]
\caption{Lowest non-dimensional vibration frequencies for a CCCC square plate
and $t=0.1$}
\label{TAB:V1}
\begin{center}
{\small\begin{tabular}{ccccccccc}
\hline
Mesh & Mode & $N=32$ & $N=64$ & $N=128$ & Order & Extrap. & \cite{HH}&\cite{DHELRS} \\
\hline
              & $\o_{11}$  &\h1.5938 &\h1.5918 &\h1.5912 &\h1.84 &1.5910 &\h1.591  &\h1.5910 \\
$\CT_{h}^{2}$ & $\o_{21}$  &\h3.0458 &\h3.0408 &\h3.0394 &\h1.85 &3.0389 &\h3.039  &\h3.0388 \\
              & $\o_{12}$  &\h3.0500 &\h3.0419 &\h3.0397 &\h1.90 &3.0389 &\h3.039  &\h3.0388 \\
              & $\o_{22}$  &\h4.2807 &\h4.2675 &\h4.2638 &\h1.85 &4.2624 &\h4.263  &\h4.2624 \\
\hline
              & $\o_{11}$  &\h1.5958 &\h1.5923 &\h1.5914 &\h1.84 &1.5910  &\h1.591 &\h1.5910\\
$\CT_{h}^{3}$ & $\o_{21}$  &\h3.0524 &\h3.0426 &\h3.0399 &\h1.83 &3.0388  &\h3.039 &\h3.0388\\
              & $\o_{12}$  &\h3.0570 &\h3.0438 &\h3.0402 &\h1.88 &3.0388  &\h3.039 &\h3.0388\\
              & $\o_{22}$  &\h4.2903 &\h4.2701 &\h4.2645 &\h1.84 &4.2623  &\h4.263 &\h4.2624\\
\hline
              & $\o_{11}$  &\h1.5961  &\h1.5923 &\h1.5914 &\h1.97 &1.5910 &\h1.591  &\h1.5910\\
$\CT_{h}^{4}$ & $\o_{21}$  &\h3.0526  &\h3.0424 &\h3.0398 &\h1.98 &3.0389 &\h3.039  &\h3.0388\\
              & $\o_{12}$  &\h3.0526  &\h3.0424 &\h3.0398 &\h1.98 &3.0389 &\h3.039  &\h3.0388\\
              & $\o_{22}$  &\h4.2914  &\h4.2699 &\h4.2644 &\h1.97 &4.2625 &\h4.263  &\h4.2624\\
\hline
              & $\o_{11}$  &\h1.5967 &\h1.5925 &\h1.5914 &\h1.98 &1.5910  &\h1.591 &\h1.5910\\
$\CT_{h}^{5}$ & $\o_{21}$  &\h3.0527 &\h3.0424 &\h3.0398 &\h1.98 &3.0389  &\h3.039 &\h3.0388\\
              & $\o_{12}$  &\h3.0573 &\h3.0435 &\h3.0401 &\h2.00 &3.0389  &\h3.039 &\h3.0388\\
              & $\o_{22}$  &\h4.2943 &\h4.2705 &\h4.2645 &\h1.98 &4.2625  &\h4.263 &\h4.2624\\
\hline
\end{tabular}}
\end{center}
\end{table}

\begin{table}[H]
\caption{Lowest non-dimensional vibration frequencies for a CCCC square plate
and $t=0.01$}
\label{TAB:V2}
\begin{center}
{\small\begin{tabular}{ccccccccc}
\hline
Mesh & Mode & $N=32$ & $N=64$ & $N=128$ & Order & Extrap. & \cite{HH}&\cite{DR} \\
\hline
              & $\o_{11}$  &\h0.1757 &\h0.1755 &\h0.1754 &\h1.89 &0.1754 &\h0.1754  &\h0.1754 \\
$\CT_{h}^{2}$ & $\o_{21}$  &\h0.3582 &\h0.3576 &\h0.3574 &\h1.87 &0.3574 &\h0.3574  &\h0.3576 \\
              & $\o_{12}$  &\h0.3587 &\h0.3577 &\h0.3575 &\h1.91 &0.3574 &\h0.3574  &\h0.3576 \\
              & $\o_{22}$  &\h0.5289 &\h0.5272 &\h0.5267 &\h1.86 &0.5265 &\h0.5264  &\h0.5274 \\
\hline
              & $\o_{11}$  &\h0.1759 &\h0.1755 &\h0.1754 &\h1.87 &0.1754 &\h0.1754  &\h0.1754 \\
$\CT_{h}^{3}$ & $\o_{21}$  &\h0.3590 &\h0.3578 &\h0.3575 &\h1.84 &0.3574 &\h0.3574  &\h0.3576 \\
              & $\o_{12}$  &\h0.3596 &\h0.3580 &\h0.3575 &\h1.87 &0.3574 &\h0.3574  &\h0.3576 \\
              & $\o_{22}$  &\h0.5304 &\h0.5276 &\h0.5268 &\h1.83 &0.5265 &\h0.5264  &\h0.5274 \\
\hline
              & $\o_{11}$  &\h0.1759 &\h0.1755 &\h0.1754 &\h1.99 &0.1754 &\h0.1754  &\h0.1754 \\
$\CT_{h}^{4}$ & $\o_{21}$  &\h0.3593 &\h0.3579 &\h0.3575 &\h2.00 &0.3574 &\h0.3574  &\h0.3576 \\
              & $\o_{12}$  &\h0.3593 &\h0.3579 &\h0.3575 &\h2.00 &0.3574 &\h0.3574  &\h0.3576 \\
              & $\o_{22}$  &\h0.5306 &\h0.5275 &\h0.5268 &\h1.99 &0.5265 &\h0.5264  &\h0.5274 \\
\hline
              & $\o_{11}$  &\h0.1762 &\h0.1756 &\h0.1754 &\h2.21 &0.1754 &\h0.1754  &\h0.1754 \\
$\CT_{h}^{5}$ & $\o_{21}$  &\h0.3597 &\h0.3579 &\h0.3575 &\h2.10 &0.3574 &\h0.3574  &\h0.3576 \\
              & $\o_{12}$  &\h0.3613 &\h0.3582 &\h0.3576 &\h2.33 &0.3574 &\h0.3574  &\h0.3576 \\
              & $\o_{22}$  &\h0.5323 &\h0.5278 &\h0.5268 &\h2.16 &0.5265 &\h0.5264  &\h0.5274 \\
\hline
\end{tabular}}
\end{center}
\end{table}

It can be seen from Tables~\ref{TAB:V1} and \ref{TAB:V2} that our method converges with a quadratic order.

Table~\ref{TAB:lockingvibra} shows the four lowest vibration frequencies computed by
Method~\ref{disc-eig} with successively refined meshes of each type for a clamped plate with
$t$=1.0e-5. The table includes orders of convegence,
as well as accurate values extrapolated by means of a least-squares fitting.
In every case, we have used a Poisson ratio $\nu=0.3$ and a correction
factor $k= 0.8601$. The reported non-dimensional frequencies are independent
of the remaining geometrical and physical parameters, except for the thickness-to-span ratio.

\begin{table}[H]
\caption{Lowest non-dimensional vibration frequencies for a CCCC square plate
and $t$=1.0e-5}
\label{TAB:lockingvibra}
\begin{center}
{\small\begin{tabular}{ccccccccc}
\hline
Mesh & Mode & $N=8$ & $N=16$ & $N=32$ & Order & Extrap. \\
\hline
              & $\o_{11}$  &\h0.7653e-1 &\h0.1289e-1 &\h0.1987e-2 &\h2.54 &-0.2995e-3  \\
$\CT_{h}^{1}$ & $\o_{21}$  &\h0.1972e-0 &\h0.3228e-1 &\h0.4942e-2 &\h2.59 &-0.5396e-3  \\
              & $\o_{12}$  &\h0.2230e-0 &\h0.4060e-1 &\h0.5085e-2 &\h2.36 &-0.3519e-2  \\
              & $\o_{22}$  &\h0.3915e-0 &\h0.6093e-1 &\h0.9893e-2 &\h2.70 & 0.7086e-3  \\
\hline
              & $\o_{11}$  &\h0.1848e-3 &\h0.1778e-3 &\h0.1761e-3 &\h2.04 &0.1756e-3  \\
$\CT_{h}^{4}$ & $\o_{21}$  &\h0.3927e-3 &\h0.3661e-3 &\h0.3601e-3 &\h2.15 &0.3583e-3  \\
              & $\o_{12}$  &\h0.3927e-3 &\h0.3661e-3 &\h0.3601e-3 &\h2.15 &0.3583e-3  \\
              & $\o_{22}$  &\h0.5983e-3 &\h0.5446e-3 &\h0.5321e-3 &\h2.11 &0.5284e-3  \\
\hline
              & $\o_{11}$  &\h0.1772e-3 &\h0.1760e-3 &\h0.1757e-3 &\h2.10 &0.1756e-3  \\
$\CT_{h}^{6}$ & $\o_{21}$  &\h0.3634e-3 &\h0.3592e-3 &\h0.3584e-3 &\h2.33 &0.3582e-3  \\
              & $\o_{12}$  &\h0.3650e-3 &\h0.3594e-3 &\h0.3584e-3 &\h2.53 &0.3582e-3  \\
              & $\o_{22}$  &\h0.5440e-3 &\h0.5312e-3 &\h0.5286e-3 &\h2.30 &0.5279e-3 \\
\hline
              & $\o_{11}$  &\h0.1779e-3 &\h0.1761e-3 &\h0.1757e-3 &\h1.94 &0.1755e-3  \\
$\CT_{h}^{7}$ & $\o_{21}$  &\h0.3654e-3 &\h0.3599e-3 &\h0.3585e-3 &\h1.95 &0.3580e-3  \\
              & $\o_{12}$  &\h0.3679e-3 &\h0.3600e-3 &\h0.3585e-3 &\h2.40 &0.3582e-3  \\
              & $\o_{22}$  &\h0.5489e-3 &\h0.5325e-3 &\h0.5289e-3 &\h2.18 &0.5278e-3  \\
\hline
\end{tabular}}
\end{center}
\end{table}

It can be seen from Table~\ref{TAB:lockingvibra} that as for the
source problem, our Method~\ref{disc-eig} lead to wrong result for
triangular meshes $\CT_{h}^{1}$ when the thickness of the plate is
small, see Remark~\ref{rem:3}. For any other family of meshes the
method is locking free and converges with a quadratic order.

Table~\ref{TAB:V3} shows the four lowest vibration frequencies computed by Method~\ref{disc-eig}
with successively refined meshes of each type for a simply supported plate with thickness
$t=0.01$. The table includes orders of convegence,
as well as accurate values extrapolated by means of a least-squares fitting.
Furthermore, the last two columns show the results reported in \cite{HH,DR}.
In every case, we have used a Poisson ratio $\nu=0.3$ and a correction
factor $k= 0.8333$. The reported non-dimensional frequencies are independent
of the remaining geometrical and physical parameters, except for the thickness-to-span ratio.

\begin{table}[H]
\caption{Lowest non-dimensional vibration frequencies for a SSSS square plate
and $t=0.01$}
\label{TAB:V3}
\begin{center}
{\small\begin{tabular}{ccccccccc}
\hline
Mesh & Mode & $N=16$ & $N=32$ & $N=64$ & Order & Extrap. & \cite{HH}&\cite{DR} \\
\hline
              & $\o_{11}$  &\h0.0966 &\h0.0963 &\h0.0963 &\h2.04 &0.0963 &\h0.0963  &\h0.0963 \\
$\CT_{h}^{2}$ & $\o_{21}$  &\h0.2416 &\h0.2408 &\h0.2406 &\h2.07 &0.2406 &\h0.2406  &\h0.2406 \\
              & $\o_{12}$  &\h0.2425 &\h0.2411 &\h0.2407 &\h2.02 &0.2406 &\h0.2406  &\h0.2406 \\
              & $\o_{22}$  &\h0.3889 &\h0.3858 &\h0.3850 &\h2.00 &0.3847 &\h0.3847  &\h0.3848 \\
\hline
              & $\o_{11}$  &\h0.0967 &\h0.0964 &\h0.0963 &\h1.96 &0.0963 &\h0.0963  &\h0.0963 \\
$\CT_{h}^{3}$ & $\o_{21}$  &\h0.2424 &\h0.2410 &\h0.2407 &\h1.91 &0.2406 &\h0.2406  &\h0.2406 \\
              & $\o_{12}$  &\h0.2434 &\h0.2413 &\h0.2408 &\h1.93 &0.2406 &\h0.2406  &\h0.2406 \\
              & $\o_{22}$  &\h0.3914 &\h0.3865 &\h0.3852 &\h1.92 &0.3847 &\h0.3847  &\h0.3848 \\
\hline
              & $\o_{11}$  &\h0.0966 &\h0.0964 &\h0.0963 &\h2.00 &0.0963 &\h0.0963  &\h0.0963 \\
$\CT_{h}^{4}$ & $\o_{21}$  &\h0.2426 &\h0.2411 &\h0.2407 &\h2.02 &0.2406 &\h0.2406  &\h0.2406 \\
              & $\o_{12}$  &\h0.2426 &\h0.2411 &\h0.2407 &\h2.02 &0.2406 &\h0.2406  &\h0.2406 \\
              & $\o_{22}$  &\h0.3898 &\h0.3860 &\h0.3850 &\h2.01 &0.3847 &\h0.3847  &\h0.3848\\
\hline
              & $\o_{11}$  &\h0.0967 &\h0.0964 &\h0.0963 &\h2.01 &0.0963 &\h0.0963  &\h0.0963 \\
$\CT_{h}^{5}$ & $\o_{21}$  &\h0.2429 &\h0.2411 &\h0.2407 &\h2.07 &0.2406 &\h0.2406  &\h0.2406 \\
              & $\o_{12}$  &\h0.2441 &\h0.2414 &\h0.2408 &\h2.09 &0.2406 &\h0.2406  &\h0.2406 \\
              & $\o_{22}$  &\h0.3910 &\h0.3863 &\h0.3851 &\h2.01 &0.3847 &\h0.3847  &\h0.3848 \\
\hline
              & $\o_{11}$  &\h0.0964 &\h0.0963 &\h0.0963 &\h2.82 &0.0963 &\h0.0963  &\h0.0963 \\
$\CT_{h}^{6}$ & $\o_{21}$  &\h0.2411 &\h0.2406 &\h0.2406 &\h3.79 &0.2406 &\h0.2406  &\h0.2406 \\
              & $\o_{12}$  &\h0.2417 &\h0.2407 &\h0.2406 &\h3.82 &0.2406 &\h0.2406  &\h0.2406 \\
              & $\o_{22}$  &\h0.3869 &\h0.3850 &\h0.3848 &\h3.40 &0.3848 &\h0.3847  &\h0.3848 \\
\hline
              & $\o_{11}$  &\h0.0965 &\h0.0963 &\h0.0963 &\h2.53 &0.0963 &\h0.0963  &\h0.0963 \\
$\CT_{h}^{7}$ & $\o_{21}$  &\h0.2416 &\h0.2408 &\h0.2406 &\h 2.35&0.2406 &\h0.2406  &\h0.2406 \\
              & $\o_{12}$  &\h0.2427 &\h0.2408 &\h0.2407 &\h3.48 &0.2406 &\h0.2406  &\h0.2406 \\
              & $\o_{22}$  &\h0.3889 &\h0.3854 &\h0.3848 &\h2.61 &0.3847 &\h0.3847  &\h0.3848 \\
\hline
\end{tabular}}
\end{center}
\end{table}

Table~\ref{TAB:V4} shows the four lowest vibration frequencies computed by Method~\ref{disc-eig}
with successively refined meshes of each type for a plate with a free edge
(with three clamped edges and the fourth free) with thickness
$t=0.01$. The table includes orders of convegence,
as well as accurate values extrapolated by means of a least-squares fitting.
Furthermore, the last two columns show the results reported in \cite{HH,DR}.
In every case, we have used a Poisson ratio $\nu=0.3$ and a correction
factor $k= 0.8601$. The reported non-dimensional frequencies are independent
of the remaining geometrical and physical parameters, except for the thickness-to-span ratio.

\begin{table}[H]
\caption{Lowest non-dimensional vibration frequencies for a CCCF square plate
and $t=0.01$}
\label{TAB:V4}
\begin{center}
{\small\begin{tabular}{ccccccccc}
\hline
Mesh & Mode & $N=32$ & $N=64$ & $N=128$ & Order & Extrap. & \cite{HH}&\cite{DR} \\
\hline
              & $\o_{11}$  &\h0.1215 &\h0.1179 &\h0.1169 &\h1.98 &0.1166 &\h0.1166  &\h0.1171  \\
$\CT_{h}^{4}$ & $\o_{21}$  &\h0.2030 &\h0.1970 &\h0.1954 &\h1.97 &0.1949 &\h0.1949  &\h0.1951  \\
              & $\o_{12}$  &\h0.3358 &\h0.3144 &\h0.3096 &\h2.14 &0.3081 &\h0.3083  &\h0.3093  \\
              & $\o_{22}$  &\h0.3884 &\h0.3773 &\h0.3745 &\h1.97 &0.3735 &\h0.3736  &\h0.3740  \\
\hline
              & $\o_{11}$  &\h0.1199 &\h0.1173 &\h0.1167 &\h2.18 &0.1166 &\h0.1166  &\h0.1171 \\
$\CT_{h}^{5}$ & $\o_{21}$  &\h0.1986 &\h0.1958 &\h0.1951 &\h2.00 &0.1948 &\h0.1949  &\h0.1951 \\
              & $\o_{12}$  &\h0.3264 &\h0.3117 &\h0.3086 &\h2.25 &0.3078 &\h0.3083  &\h0.3093 \\
              & $\o_{22}$  &\h0.3791 &\h0.3749 &\h0.3738 &\h1.92 &0.3734 &\h0.3736  &\h0.3740 \\
\hline
              & $\o_{11}$  &\h0.1177 &\h0.1169 &\h0.1167 &\h2.00 &0.1166 &\h0.1166  &\h0.1171 \\
$\CT_{h}^{6}$ & $\o_{21}$  &\h0.1967 &\h0.1953 &\h0.1949 &\h2.13 &0.1948 &\h0.1949  &\h0.1951 \\
              & $\o_{12}$  &\h0.3134 &\h0.3090 &\h0.3081 &\h2.22 &0.3079 &\h0.3083  &\h0.3093 \\
              & $\o_{22}$  &\h0.3753 &\h0.3738 &\h0.3735 &\h2.30 &0.3734 &\h0.3736  &\h0.3740 \\
\hline
              & $\o_{11}$  &\h0.1180 &\h0.1169 &\h0.1167 &\h1.90 &0.1166 &\h0.1166  &\h0.1171 \\
$\CT_{h}^{7}$ & $\o_{21}$  &\h0.1974 &\h0.1954 &\h0.1950 &\h2.15 &0.1948 &\h0.1949  &\h0.1951 \\
              & $\o_{12}$  &\h0.3151 &\h0.3095 &\h0.3082 &\h2.08 &0.3078 &\h0.3083  &\h0.3093 \\
              & $\o_{22}$  &\h0.3772 &\h0.3743 &\h0.3736 &\h2.22 &0.3734 &\h0.3736  &\h0.3740 \\
\hline
\end{tabular}}
\end{center}
\end{table}

It can be seen from Tables~\ref{TAB:V3} and \ref{TAB:V4} that our method converges with a quadratic order.

\subsection{Buckling of plates}

The effectiveness of the MDF method for buckling analysis are demostrated
by examples with different thickness, boundary conditions and different in-plane
compressive stress $\bsi$.

We have computed approximations of the buckling coefficients $\l^{bc}=\l^{bp}t^{2}$ being the smallest
(the critical load) by which the chosen in-plane
compressive stress $\bsi$ must be multiplied by in order to cause buckling.
In order to compare our results with
those in \cite{KXWL,LC,NLTN}, a non-dimensional buckling intensity is defined as:
$$K:=\frac{\l_{h}^{bc}L}{\pi^2 D},$$
here $\l_{h}^{bc}=\l_{h}^{bp}t^{2}$ are the computed buckling coefficients, $L$ is the plate side length
and $D$ is the flexural rigidity defined as $D=Et^3/[12(1-\nu^2)]$.

\subsubsection{Uniformly compressed plate}\label{ss:unif-plate}

In this couple of tests, we use $\bsi={\bf I}$, corresponding to a
uniformly compressed plate (in the $x$, $y$ directions).

First, we consider a simply supported plate, since analytical solutions
are available (see~\cite{R2}) for that case.
In Table~\ref{TAB:B1}, we report the four lowest non-dimensional
buckling intensities $K_{1},\ldots,K_{4}$, for the thickness $t=0.01$, and $L=1$
computed by Method~\ref{disc-stress} with four different family of meshes.
The table includes computed orders of convergence,
as well as more accurate values extrapolated
by means of a least-squares procedure. Furthermore, the last column reports
the exact buckling intensities. In this case, we have used a Poisson ratio $\nu=0.3$ and a correction
factor $k= 5/6$.

\begin{table}[H]
\caption{Lowest non-dimensional buckling intensities $K_{1},\ldots,K_{4}$ for a SSSS square plate and $t=0.01$}
\label{TAB:B1}
\begin{center}
{\small\begin{tabular}{cccccccc}
\hline
Mesh & $K$ & $N=16$ & $N=32$ & $N=64$ & Order & Extrap. & Exact\\
\hline
              & $K_{1}$  &\h2.0315 &\h2.0075 &\h2.0011 &\h1.90 &1.9987 &\h1.9989  \\
$\CT_{h}^{3}$ & $K_{2}$  &\h5.1607 &\h5.0381 &\h5.0046 &\h1.87 &4.9920 &\h4.9930  \\
              & $K_{3}$  &\h5.2262 &\h5.0541 &\h5.0086 &\h1.92 &4.9922 &\h4.9930 \\
              & $K_{4}$  &\h8.5170 &\h8.1249 &\h8.0186 &\h1.88 &7.9788 &\h7.9820 \\
\hline
              & $K_{1}$  &\h2.0381 &\h2.0086 &\h2.0013 &\h2.01 &1.9989 &\h1.9989  \\
$\CT_{h}^{4}$ & $K_{2}$  &\h5.2116 &\h5.0465 &\h5.0063 &\h2.04 &4.9934 &\h4.9930  \\
              & $K_{3}$  &\h5.2116 &\h5.0465 &\h5.0063 &\h2.04 &4.9934 &\h4.9930  \\
              & $K_{4}$  &\h8.6292 &\h8.1388 &\h8.0209 &\h2.06 &7.9839 &\h7.9820 \\
\hline
              & $K_{1}$  &\h2.0412 &\h2.0093 &\h2.0015 &\h2.02 &1.9989 &\h1.9989  \\
$\CT_{h}^{5}$ & $K_{2}$  &\h5.2242 &\h5.0485 &\h5.0063 &\h2.06 &4.9931 &\h4.9930  \\
              & $K_{3}$  &\h5.2929 &\h5.0641 &\h5.0099 &\h2.08 &4.9932 &\h4.9930  \\
              & $K_{4}$  &\h8.6788 &\h8.1504 &\h8.0234 &\h2.06 &7.9836 &\h7.9820 \\
\hline
              & $K_{1}$  &\h2.0347 &\h2.0068 &\h2.0010 &\h2.27 &1.9995 &\h1.9989  \\
$\CT_{h}^{7}$ & $K_{2}$  &\h5.1962 &\h5.0424 &\h5.0053 &\h2.05 &4.9935 &\h4.9930  \\
              & $K_{3}$  &\h5.2429 &\h5.0451 &\h5.0063 &\h2.35 &4.9968 &\h4.9930   \\
              & $K_{4}$  &\h8.5851 &\h8.1192 &\h8.0158 &\h2.17 &7.9862 &\h7.9820 \\
\hline
\end{tabular}}
\end{center}
\end{table}

It can be seen from Table~\ref{TAB:B1} that our method converges to the exact values
with a quadratic order.

As a second test, we present the results for the lowest non-dimensional
buckling intensity $K_{1}$ for a clamped plate with varying
thickness $t$, in order to assess the stability of the Method~\ref{disc-stress}
when $t$ goes to zero. It is well known that $K_{1}$ converges
to the non-dimensional buckling intensity of an identical Kirchhoff-Love
uniformly compressed clamped plate.

In Table~\ref{TAB:B2}, we report the lowest non-dimensional buckling intensity
$K_{1}$ of a uniformly compressed clamped plate with varying thickness $t$ and $L=1$.
We have used five different family of meshes.
The table includes computed orders of convergence,
as well as more accurate values extrapolated
by means of a least-squares procedure.
In the last row of each family of meshes
we report the limit values as $t$ goes to zero obtained by extrapolation.
In this case, we have used a Poisson ratio $\nu=0.25$ and a correction
factor $k= 5/6$.

\begin{table}[H]
\caption{Lowest non-dimensional buckling intensity $K_{1}$ of a clamped plate with varying
thickness.}
\label{TAB:B2}
\begin{center}
{\small\begin{tabular}{|c||c||c|c|c|c|c|}
\hline
Mesh & $t$ & $N=16$ & $N=32$ & $N=64$ & Order & Extrap. \\
\hline
              & 0.1     &\h4.9031 &\h4.6658 &\h4.6099 &\h2.09 &4.5929  \\
$\CT_{h}^{1}$ & 0.01    &\h7.1314 &\h5.5307 &\h5.3275 &\h2.98 &5.2982 \\
              & 0.001   &\h9.9817e+1 &\h9.2546 &\h5.6315 &\h4.00  & 4.3035 \\
              & 0.0001  &\h9.2723e+3 &\h2.6878e+2 &\h1.2923e+1 &\h4.00 &-0.1678 \\
\hline
              & 0 (extrap.)  &\h --&\h --&\h-- &\h-- &-- \\
\hline
\hline
              & 0.1     &\h4.6316 &\h4.6015 &\h4.5937 &\h1.93 &4.5909  \\
$\CT_{h}^{2}$ & 0.01    &\h5.3423 &\h5.3072 &\h5.2981 &\h1.96 &5.2950 \\
              & 0.001   &\h5.3509 &\h5.3157 &\h5.3066 &\h1.96 &5.3035  \\
              & 0.0001  &\h5.3510 &\h5.3158 &\h5.3067 &\h1.96 &5.3035 \\
\hline
              & 0 (extrap.)  &\h5.3510 &\h5.3158 &\h5.3067 &\h1.96 &5.3036 \\
\hline
\hline
              & 0.1     &\h4.6476 &\h4.6058 &\h4.5948 &\h1.92 &4.5908  \\
$\CT_{h}^{3}$ & 0.01    &\h5.3611 &\h5.3121 &\h5.2994 &\h1.94 &5.2949  \\
              & 0.001   &\h5.3697 &\h5.3207 &\h5.3079 &\h1.94 &5.3034  \\
              & 0.0001  &\h5.3698 &\h5.3208 &\h5.3080 &\h1.94 &5.3035 \\
\hline
              & 0 (extrap.)  &\h5.3698 &\h5.3208 &\h5.3080 &\h1.94 &5.3035 \\
\hline
\hline
              & 0.1     &\h4.6441 &\h4.6043 &\h4.5943 &\h2.00 & 4.5910  \\
$\CT_{h}^{4}$ & 0.01    &\h5.3564 &\h5.3103 &\h5.2989 &\h2.01 &  5.2951\\
              & 0.001   &\h5.3649 &\h5.3188 &\h5.3074 &\h2.01 & 5.3036 \\
              & 0.0001  &\h5.3649 &\h5.3189 &\h5.3074 &\h2.01 & 5.3037\\
\hline
              & 0 (extrap.)  &\h5.3650 &\h5.3189 &\h5.3074 &\h2.01 &5.3037 \\
\hline
\hline
              & 0.1     &\h4.6487 &\h4.6054 &\h4.5946 &\h2.00 &4.5909  \\
$\CT_{h}^{5}$ & 0.01    &\h5.3702 &\h5.3128 &\h5.2993 &\h2.09 &5.2952  \\
              & 0.001   &\h5.3863 &\h5.3240 &\h5.3085 &\h2.01 &5.3034  \\
              & 0.0001  &\h5.3866 &\h5.3242 &\h5.3088 &\h2.01 &5.3036 \\
\hline
              & 0 (extrap.)  &\h5.3867 &\h5.3242 &\h5.3087 &\h2.01 &5.3036 \\
\hline
\end{tabular}}
\end{center}
\end{table}

Additionally, we have also computed the lowest buckling intensity
of a Kirchhoff-Love plate by using the finite element method analyzed in \cite{MR}.

In Table~\ref{TAB:B3}, we report the lowest non-dimensional buckling intensity
of a uniformly compressed clamped plate with $L=1$.
In this case we considered a Poisson ratio $\nu=0.25$.

\begin{table}[H]
\caption{Lowest non-dimensional buckling intensity of a uniformly compressed clamped
thin plate (Kirchhoff-Love model) computed with the method from~\cite{MR}.}
\label{TAB:B3}
\begin{center}
{\small\begin{tabular}{ccccccc}
\hline
 Method & $N=24$ & $N=36$ & $N=48$ & $N=60$ & Order & Extrapolated \\
\hline
        \cite{MR}  &\h5.3051 &\h5.3042 &\h5.3039 &\h5.3038 & 2.61 &\h5.3037 \\
\hline
\end{tabular}}
\end{center}
\end{table}

It is clear from the results of Tables~\ref{TAB:B2} and \ref{TAB:B3},
that our Method~\ref{disc-stress} lead to wrong result for triangular
meshes $\CT_{h}^{1}$ when the thickness of the plate is small, see Remark~\ref{rem:3}.
For all the other family of meshes the method is locking free and do not
deteriorate as the plate thickness become smaller.

\subsubsection{Plate uniformly compressed in one direction}\label{unif-plate1direc}

In this couple of tests, we use
$$\bsi=\left[\begin{array}{cc}1&0\\0&0\end{array}\right],$$
corresponding to a plate subjected to uniaxial compression
(in the $x$ direction).
We consider different boundary conditions.

In Tables~\ref{TAB:B4} and \ref{TAB:B5}, we report the lowest non-dimensional buckling intensity
$K_{1}$, for a clamped and simply supported plate, respectively, with thickness $t=0.1$, and $L=1$
computed by Method~\ref{disc-stress} with different family of meshes.
The table includes computed orders of convergence,
as well as more accurate values extrapolated
by means of a least-squares procedure.
Furthermore, the last two columns show the results reported in \cite{KXWL,LC}.
In these cases, we have used a Poisson ratio $\nu=0.3$ and a correction
factor $k= 5/6$.

\begin{table}[H]
\caption{Lowest non-dimensional buckling intensity $K_{1}$ for a CCCC square plate and $t=0.1$}
\label{TAB:B4}
\begin{center}
{\small\begin{tabular}{ccccccccc}
\hline
Mesh & $K$ & $N=32$ & $N=64$ & $N=128$ & Order & Extrap. & \cite{KXWL}&\cite{LC} \\
\hline
$\CT_{h}^{2}$ & $K_{1}$  &\h8.3849 &\h8.3157 &\h8.2978 &\h1.95 &8.2915 &\h8.2917 &\h8.2931  \\
\hline
$\CT_{h}^{3}$ & $K_{1}$  &\h8.4222 &\h8.3260 &\h8.3004 &\h1.91 &8.2912 &\h8.2917  &\h8.2931 \\
\hline
$\CT_{h}^{4}$ & $K_{1}$  &\h8.3987 &\h8.3185 &\h8.2984 &\h2.00 &8.2917 &\h8.2917  &\h8.2931 \\
\hline
$\CT_{h}^{5}$ & $K_{1}$  &\h8.4273 &\h8.3255 &\h8.3001 &\h2.00 &8.2916 &\h8.2917  &\h8.2931 \\
\hline
$\CT_{h}^{6}$ & $K_{1}$  &\h8.3663 &\h8.3110 &\h8.2963 &\h1.91 &8.2910 &\h8.2917  &\h8.2931 \\
\hline
$\CT_{h}^{7}$ & $K_{1}$  &\h8.3715 &\h8.3121 &\h8.2965 &\h1.93 &8.2909 &\h8.2917  &\h8.2931 \\
\hline
\end{tabular}}
\end{center}
\end{table}

\begin{table}[H]
\caption{Lowest non-dimensional buckling intensity $K_{1}$ for a SSSS square plate and $t=0.1$}
\label{TAB:B5}
\begin{center}
{\small\begin{tabular}{ccccccccc}
\hline
Mesh & K & $N=32$ & $N=64$ & $N=128$ & Order & Extrap. & \cite{KXWL}&\cite{LC} \\
\hline
$\CT_{h}^{2}$ & $K_{1}$  &\h3.7993  &\h3.7897 &\h3.7873 &\h1.97 &3.7864 &\h3.7865  &\h3.7873 \\
\hline
$\CT_{h}^{3}$ & $K_{1}$  &\h3.8029 &\h3.7907 &\h3.7875  &\h1.96 &3.7864 &\h3.7865  &\h3.7873 \\
\hline
$\CT_{h}^{4}$ & $K_{1}$  &\h3.8049 &\h3.7911 &\h3.7876 &\h2.00 &3.7864 &\h3.7865  &\h3.7873 \\
\hline
$\CT_{h}^{5}$ & $K_{1}$  &\h3.8062 &\h3.7913 &\h3.7877 &\h2.01 &3.7865 &\h3.7865  &\h3.7873 \\
\hline
$\CT_{h}^{6}$ & $K_{1}$  &\h3.7975 &\h3.7892 &\h3.7871 &\h1.98 &3.7864 &\h3.7865  &\h3.7873 \\
\hline
$\CT_{h}^{7}$ & $K_{1}$  &\h3.8011 &\h3.7903 &\h3.7874 &\h1.90 &3.7863 &\h3.7865  &\h3.7873 \\
\hline
\end{tabular}}
\end{center}
\end{table}

It can be seen from Tables~\ref{TAB:B4} and \ref{TAB:B5} that our
method converges with a quadratic order.

\subsubsection{Shear loaded plate}\label{shear-plate}

In this test, we use
$$\bsi=\left[\begin{array}{cc}0&1\\1&0\end{array}\right],$$
corresponding to a plate subjected to shear load.
We consider different boundary conditions.

In Table~\ref{TAB:B6}, we report the lowest non-dimensional buckling intensity
$K_{1}$, for a simply supported plate with thickness $t=0.01$, and $L=1$
computed by Method~\ref{disc-stress} with different family of meshes.
The table includes computed orders of convergence,
as well as more accurate values extrapolated
by means of a least-squares procedure.
Furthermore, the last two columns show the results reported in \cite{NLTN}.
In these cases, we have used a Poisson ratio $\nu=0.3$ and a correction
factor $k= 5/6$.

\begin{table}[H]
\caption{Lowest non-dimensional buckling intensity $K_{1}$ for a SSSS square plate and $t=0.01$}
\label{TAB:B6}
\begin{center}
{\small\begin{tabular}{cccccccc}
\hline
Mesh & $K$ & $N=32$ & $N=64$ & $N=128$ & Order & Extrap. & \cite{NLTN} \\
\hline
$\CT_{h}^{2}$ & $K_{1}$  &\h9.4832 &\h9.3514 &\h9.3180 &\h1.98 &9.3067 &\h9.2830 \\
\hline
$\CT_{h}^{3}$ & $K_{1}$  &\h9.3848 &\h9.3270 &\h9.3119 &\h1.94 &9.3066 &\h9.2830 \\
\hline
$\CT_{h}^{4}$ & $K_{1}$  &\h9.4602 &\h9.3450 &\h9.3164 &\h2.01 &9.3069 &\h9.2830 \\
\hline
$\CT_{h}^{5}$ & $K_{1}$  &\h9.4759 &\h9.3483 &\h9.3170 &\h2.03 &9.3069 &\h9.2830 \\
\hline
$\CT_{h}^{6}$ & $K_{1}$  &\h9.4196 &\h9.3346 &\h9.3134 &\h2.01 &9.3064 &\h9.2830\\
\hline
$\CT_{h}^{7}$ & $K_{1}$  &\h9.4640 &\h9.3460 &\h9.3163 &\h1.99 &9.3063 &\h9.2830 \\
\hline
\end{tabular}}
\end{center}
\end{table}

It can be seen from Table~\ref{TAB:B6} that our
method converges with a quadratic order.

\section{Conclusions}

We assessed numerically the actual performance of the method proposed
in \cite{BM}, extending it also to free vibration and buckling problems of plates.
We tested different families of mimetic meshes, different values of the
relative thickness and various boundary conditions. In all the three
types of problems considered (source problem, free vibration, buckling)
the method was shown to be locking free and to converge with an optimal
rate both in discrete $L^\infty$ and $H^1$ norms for meshes
made with elements with 4 or more edges. In some occasions, a
super convergence rate was noticed. Moreover, differently from standard
quadrilateral finite elements, the method shows a robust behavior also
for uniformly distorted families of meshes such as those in
Figure~\ref{FIG:UNIF_MESH5}. We thus conclude that the proposed
method is very reliable for Reissner-Mindlin plate computations.

\begin{acknowledgements}
The third author was partially supported by CONICYT-Chile through FONDECYT project No. 11100180, and
by Centro de Investigaci\'on en Ingenier\'\i a Matem\'atica (CI$^2$MA), Universidad de Concepci\'on, Chile.
\end{acknowledgements}


\end{document}